\documentclass[a4paper,12pt]{article}
\usepackage[utf8]{inputenc}
\usepackage[english]{babel}

\usepackage{adjustbox}
\usepackage{amssymb}
\usepackage{amsmath}
\usepackage{amssymb}
\usepackage[affil-it]{authblk}
\usepackage[sorting=none,url=false,isbn=false,style=ieee]{biblatex}
\usepackage{booktabs}
\usepackage{csquotes}
\usepackage{changepage}
\usepackage{float}
\usepackage{geometry}
\usepackage[acronyms,nonumberlist]{glossaries-extra}
\usepackage{hyperref}
\usepackage[capitalize]{cleveref}
\usepackage{mathtools}
\usepackage{multicol}
\usepackage{multirow}
\usepackage{physics}
\usepackage{rotating}
\usepackage{siunitx}
\usepackage{subcaption}
\usepackage[dvipsnames]{xcolor}
\usepackage{tikz}

\hypersetup{colorlinks,allcolors=cyan}
\setabbreviationstyle[acronym]{long-short}
\newacronym{ads}{ADS}{automatic domain splitting}
\newacronym{aop}{AOP}{argument of periapsis}

\newacronym{ci}{CI}{confidence interval}
\newacronym{cnes}{CNES}{Centre National d'\'Etudes Spatiales}
\newacronym{cut}{CUT}{conjugate unscented transform}

\newacronym{da}{DA}{differential algebra}
\newacronym{de}{DE}{Developmental Ephemerides}

\newacronym{fp}{FP}{floating point}

\newacronym{gn}{GN}{Gauss-Newton}
\newacronym{gto}{GTO}{geostationary transfer orbit}

\newacronym{hf}{HF}{high-fidelity}

\newacronym{ic}{IC}{initial condition}
\newacronym{iod}{IOD}{initial orbit determination}

\newacronym{jpl}{JPL}{Jet Propulsion Laboratory}

\newacronym{kf}{KF}{Kalman filter}

\newacronym{lc}{LinCov}{linear covariance}
\newacronym{lf}{LF}{low-fidelity}
\newacronym{lm}{LM}{Levenberg–Marquardt}
\newacronym{loads}{LOADS}{low-order automatic domain splitting}
\newacronym{los}{LOS}{line of sight}
\newacronym{lpp}{LPP}{linear programming problem}
\newacronym{ls}{LS}{least squares}
\newacronym{lsar}{LSAR}{least sum of absolute residuals}

\newacronym{mc}{MC}{Monte Carlo}
\newacronym{mf}{MF}{multifidelity}

\newacronym{nli}{NLI}{nonlinearity index}
\newacronym{norad}{NORAD}{North American Aerospace Defense Command}

\newacronym{od}{OD}{orbit determination}
\newacronym{ode}{ODE}{ordinary differential equation}

\newacronym{pace}{PACE}{Polynomial Algebra Computational Engine}
\newacronym{pce}{PCE}{polynomial chaos expansion}
\newacronym{pdf}{pdf}{probability density function}
\newacronym{pn}{PN}{process noise}
\newacronym{pw}{PW}{point-wise}

\newacronym{ra}{RA}{right ascension}
\newacronym{raan}{RAAN}{right ascension of the ascending node}
\newacronym{rmse}{RMSE}{root mean square error}

\newacronym{sgp}{SGP4}{Simplified General Perturbations}
\newacronym{snc}{SNC}{State Noise Compensation}
\newacronym{so}{SO}{space object}
\newacronym{srp}{SRP}{Solar radiation pressure}
\newacronym{ssa}{SSA}{space situational awareness}
\newacronym{stm}{STM}{state transition matrix}
\newacronym{stt}{STT}{State Transition Tensor}

\newacronym{tarot}{TAROT}{T\'elescope \`a Action Rapide pour les Objets Transitoires}
\newacronym{tle}{TLE}{Two-Line Element}

\newacronym{up}{UP}{uncertainty propagation}
\newacronym{ut}{UT}{unscented transform}

\newcommand{\glsentryfullinv}[1]{\glsentryshort{#1} (\glsentrylong{#1})}
\DeclareMathOperator*{\argmax}{arg\,max}

\usetikzlibrary{shapes.geometric,arrows}

\tikzstyle{startstop} = [rectangle, rounded corners, 
minimum width=3cm, 
minimum height=1cm,
text centered, 
draw=black, 
fill=red!30]

\tikzstyle{io} = [trapezium, 
trapezium stretches=true,
trapezium left angle=70, 
trapezium right angle=110, 
minimum width=3cm, 
minimum height=1cm, text centered, 
draw=black, fill=blue!30]

\tikzstyle{process} = [rectangle,
minimum width=3cm, 
minimum height=1cm, 
text centered, 
text width=5cm, 
draw=black, 
fill=orange!30]

\tikzstyle{smallprocess} = [rectangle,
minimum width=3cm, 
minimum height=1cm, 
text centered, 
text width=3cm, 
draw=black, 
fill=orange!30]

\tikzstyle{decision} = [diamond, 
minimum width=3cm, 
minimum height=1cm, 
text centered, 
draw=black, 
fill=green!30]

\tikzstyle{arrow} = [thick,->,>=stealth]

\makeatletter
\def\@maketitle{
\begin{center}
    {\Large\@title}\\[2em]
    {\small\@author}
\end{center}
}
\makeatother

\title{A Multifidelity Approach to Robust Orbit Determination}

\author[1,*]{Alberto Foss\`a}
\author[2]{Roberto Armellin}
\author[3]{Emmanuel Delande}
\author[1]{Matteo Losacco}
\author[1]{Francesco Sanfedino}

\affil[1]{Department of Aerospace Vehicles Design and Control, Institut Supérieur de l'Aéronautique et de l’Espace, 10 Avenue Edouard Belin, Toulouse, 31055, France}
\affil[2]{Te Pūnaha Ātea - Space Institute, The University of Auckland, 20 Symonds Street, Auckland, 1010, New Zealand}
\affil[3]{Space Situational Awareness Office, Centre National d'Études Spatiales, 18 Avenue Edouard Belin, Toulouse, 31400, France}
\affil[*]{Corresponding author, \href{mailto:alberto.fossa@isae-supaero.fr}{alberto.fossa@isae-supaero.fr}}

\begin{document}

\maketitle

\begin{abstract}

This paper presents an algorithm for the preprocessing of observation data aimed at improving the robustness of \glsentrylong{od} tools. Two objectives are fulfilled: obtain a refined solution to the \glsentrylong{iod} problem and detect possible outliers in the processed measurements. The uncertainty on the initial estimate is propagated forward in time and progressively reduced by exploiting sensor data available in said propagation window. \Glsentrylong{da} techniques and a novel \glsentrylong{ads} algorithm for second-order Taylor expansions are used to efficiently propagate uncertainties over time. A \glsentrylong{mf} approach is employed to minimize the computational effort while retaining the accuracy of the propagated estimate. At each observation epoch, a polynomial map is obtained by projecting the propagated states onto the observable space. Domains that do no overlap with the actual measurement are pruned thus reducing the uncertainty to be further propagated. Measurement outliers are also detected in this step. The refined estimate and retained observations are then used to improve the robustness of batch \glsentrylong{od} tools. The effectiveness of the algorithm is demonstrated for a \glsentrylong{gto} object using synthetic and real observation data from the \glsentryshort{tarot} network.

\end{abstract}


\section{Introduction}

Preserving the realism of the probabilistic representation of \glsentryfullpl{so} states over time is paramount in \glsentryfull{ssa} applications such as \glsentryfull{od} and \glsentryshortpl{so} catalogs build-up and maintenance. To accomplish these tasks, two main tools are required: an accurate and efficient \glsentryfull{up} method and an effective strategy to incorporate newly available information in the predicted estimate. This work contributes to these needs by proposing a pruning algorithm for the refinement of an initial solution, typically the result of an \glsentryfull{iod} process, and the detection of measurements outliers among the input set. The algorithm is composed of four main steps which are carried out sequentially for each available measurements: propagation of the current state estimate, projection onto the observable space, pruning of the uncertainty based on the newly available information, and merging to reduce the information to be further propagated.

Regarding the \glsentryshort{up} problem in orbital dynamics, this has been historically tackled with either linearized methods or large scale \glsentryfull{mc} simulations. Several methods for nonlinear \glsentryshort{up} have since been proposed and are conveniently classified as intrusive or non-intrusive, with the first requiring full knowledge of the underlying dynamics while the seconds treating the dynamical model as a \emph{black box}. The first category includes \glsentryfullpl{stt}~\cite{Park2006a} and \glsentryfull{da}~\cite{Berz1999,Valli2013} while the second gathers the aforementioned \glsentryshort{mc} as well as \glsentryfull{ut}~\cite{Julier2004}, \glsentryfull{cut}~\cite{Adurthi2015} and \glsentryfullpl{pce}~\cite{Jones2013}. A comprehensive survey of these methods can be found in \cite{Luo2017}. In this work, \glsentryshort{da} is used for \glsentryshort{up} in conjunction with an \glsentryfull{ads} algorithm~\cite{Wittig2015a} to control the error on the final solution manifold. \Glsentrylong{mf} methods are exploited to drastically reduce the computational effort of the overall algorithm while retaining a similar accuracy on the final solution. Their effectiveness for orbit \glsentryshort{up} was already demonstrated in~\cite{Jones2018,Fossa2022} and corroborated in this work. The \glsentryshort{da}-based \glsentryshort{up} is carried out on a \glsentrylong{lf} dynamical model and only the centers of the resulting manifold are propagated in \glsentrylong{hf} to compute a posteriori correction to the \glsentrylong{lf} solution.

The projection onto the observables space follows, also wrapped by the \glsentryshort{ads} algorithm to retain the accuracy of the resulting polynomial map, or manifold, that represents the uncertainty on the predicted measurements by patching a set of truncated Taylor expansions. Polynomial bounding techniques are then employed to estimate the range of variation of both predicted and actual measurements. These bounds are then used to discard the domains that do not match the real observation and reduce the size of the manifold of states to be further propagated. Possible measurements outliers are also identified in the observations that have no intersection with the projected set. The output is a refined solution to the \glsentryshort{iod} problem that is used as improved guess to initialize a batch \glsentryshort{od} algorithm. Concerning the last, if \glsentryfull{ls} techniques are the workhorse of current operational algorithms, estimators that minimize the $L_1$ norm of the residuals vector were shown to improve robustness in the presence of measurements outliers~\cite{Narula1982}. These \glsentryfull{lad} estimators were also proven to be faster than \glsentryshort{ls} regression on large data sets~\cite{Portnoy1997} and successfully applied to the problem of batch \glsentryshort{od}~\cite{Branham2005,Gehly2016}. The \glsentryfull{lsar} algorithm was then proposed in~\cite{Prabhu2021} as an alternative technique to compute the minimum $L_1$ norm estimates. This algorithm recasts the \glsentryshort{lad} problem into a sequence of \glsentryfullpl{lpp} which can be then solved with highly efficient simplex or interior point methods~\cite{Nocedal2006}. Both solutions are compared in this paper highlighting the benefits provided by the pruning algorithm over a direct processing of all available information. A complete solution to the batch \glsentryshort{od} problem is thus presented, with \glsentryshort{da} being the common denominator underlying each proposed tool. This approach is demonstrated of being robust to measurements outliers thanks to an effective preprocessing of the available observations and the use of the \glsentryshort{lsar} algorithm to compute the final estimate.

The paper is organized as follows. \Cref{sec:math} sets the background material for the paper. This includes the modeling of the problem in \cref{sec:modeling}, the \glsentryshort{da} framework in \cref{sec:da} and the \glsentryshort{od} tools in \cref{sec:od}. The main contribution of the current work, namely the \glsentrylong{mf} pruning scheme, is then detailed in \cref{sec:robust_od} after  the \glsentryshort{da}-based \glsentryshort{iod} algorithm is briefly outlined. Numerical applications of the proposed tool are presented in \cref{sec:application} before drawing the conclusions in \cref{sec:conclusions}.

\section{Mathematical background}\label{sec:math}

This section describes the building blocks this work is built on. These include the modeling of the \glsentrylong{pn} and the measurement process, \glsentrylong{da} techniques applied to \glsentrylong{up}, the \glsentrylong{ads} and merging algorithms used to control the accuracy of Taylor expansions, and the batch \glsentrylong{od} tools.

\subsection{Modeling of the problem}\label{sec:modeling}
The numerical algorithm used to model the \glsentrylong{pn} is introduced hereafter, followed by a description of the measurement process and sensor noise model for optical telescopes.
\subsubsection{State noise compensation}\label{sec:snc}
Taking into account the uncertainty due to neglected or mismodeled accelerations is paramount to achieve a precise estimate of the propagated orbit. Stochastic accelerations are thus included in the developed \glsentryshort{up} framework and their effects on the state covariance estimated using the \glsentryfull{snc} algorithm~\cite{Tapley2004}.

For continuous systems, the dynamics under the influence of \glsentrylong{pn} is described by
\begin{equation}
    \dot{\vb*{x}} = \vb*{f}(\vb*{x},t) + \vb*{B}(t)\vb*{w}(t)
    \label{eqn:ode_noise}
\end{equation}
with $\vb*{f}:\mathbb{R}^n\times\mathbb{R}\to\mathbb{R}^n$ the deterministic dynamics, $\vb*{w}(t)$ the $m$-dimensional \glsentrylong{pn}, and $\vb*{B}(t)$ the $n\times m$ \emph{\glsentrylong{pn} mapping matrix} that maps the \glsentrylong{pn} into state derivatives. Under the assumption of uncorrelated white noise, the first two statistical moments of $\vb*{w}(t)$ are given by
\begin{equation}
    \begin{aligned}
        \mathbb{E}[\vb*{w}(t)] &= \vb*{0}\\
        \mathbb{E}[\vb*{w}(t)\vb*{w}^T(t)] &= \vb*{Q}(t)\delta(t-\tau)
    \end{aligned}
    \label{eqn:white_noise_exp}
\end{equation}
with $\delta(t-\tau)$ the Dirac delta and $\vb*{Q}(t)$ the \glsentrylong{psd} function of $\vb*{w}(t)$~\cite{Carpenter2018}. An estimate of the true state is then maintained through a mean state $\vb*{x}(t)$ and a covariance matrix $\vb*{P}(t)$ that are computed by integrating the system of first-order \glsentryfullpl{ode}
\begin{equation}
    \begin{dcases}
        \dot{\vb*{x}}(t) &= \vb*{f}(\vb*{x},t)\\
        \dot{\vb*{P}}(t) &= \vb*{A}(t)\vb*{P}(t) + \vb*{P}(t)\vb*{A}^T(t) + \vb*{B}(t)\vb*{Q}(t)\vb*{B}^T(t)
    \end{dcases}
    \label{eqn:ode_snc}
\end{equation}
subject to \glsentrylongpl{ic} $\vb*{x}(t_0) = \vb*{x}_0$ and $\vb*{P}(t_0) = \vb*{0}$.

\subsubsection{Measurement process and sensor noise}\label{sec:meas_model}
Even though the pruning algorithm presented in \cref{sec:mf_up} is agnostic with respect to the type of measurements, simulations in \cref{sec:application} uses angular measurements from optical sensors. The corresponding measurement model is thus briefly introduced hereafter.

Consider an optical telescope identified by its geodetic coordinates $(\phi,\lambda,h)$ with $\phi\in[-\pi/2,\pi/2]$ geodetic latitude, $\lambda\in[-\pi,\pi]$ geodetic longitude and $h\in\mathbb{R}^+$ geodetic height. Denote with $\vb*{x}=[\vb*{r}^T\ \vb*{v}^T]^T$ the inertial state vector of the tracked \glsentryshort{so} and with $\vb*{r}_{obs}$ the inertial position of the telescope at the same epoch. The two positions $\vb*{r},\vb*{r}_{obs}$ are related by
\begin{equation}
    \vb*{r} = \vb*{r}_{obs}+\rho\hat{\vb*{\rho}}
    \label{eq:so_pos_from_range}
\end{equation}
with $\rho$ topocentric range and $\hat{\vb*{\rho}}$ \glsentryfull{los} unit vector. The last is computed from the topocentric right ascension and declination $(\alpha,\delta)$ as
\begin{equation}
    \hat{\vb*{\rho}} = \begin{bmatrix}
        \cos{\alpha}\cos{\delta}\\
        \sin{\alpha}\cos{\delta}\\
        \sin{\delta}
    \end{bmatrix}
    \label{eq:sph2los}
\end{equation}

If the \glsentryshort{los} vector $\rho\hat{\vb*{\rho}}=[\rho_x\ \rho_y\ \rho_z]^T$ is known instead, the topocentric range, right ascension and declination are given by
\begin{equation}
    \begin{dcases}
        \rho = \sqrt{\rho_x^2+\rho_y^2+\rho_z^2}\\
        \alpha = \tan^{-1}{\rho_y/\rho_x}\\
        \delta = \sin^{-1}{\rho_z/\rho}
    \end{dcases}
    \label{eq:cart2sph}
\end{equation}

Sensor noise is modeled by assuming uncorrelated angular measurements $(\alpha_i,\delta_i)$ corrupted by Gaussian white noise with standard deviations $(\sigma_{\alpha_i},\sigma_{\delta_i})$ for each measurement epoch $t_i$. These uncertainties are represented in the \glsentryshort{da} framework as
\begin{equation}
\begin{aligned}
    [\alpha_i] &= \alpha_i + c\sigma_{\alpha_i}\cdot\delta\alpha_i\\
    [\delta_i] &= \delta_i + c\sigma_{\delta_i}\cdot\delta\delta_i
\end{aligned}
\label{eq:ang_meas_da}
\end{equation}
where $c>0$ is the $z$-score that quantifies how many standard deviations lie within the polynomial bounds. This parameter is typically set equal to $c=3$ such that the true measurement is found within the domain of the polynomial with a $\approx 99.73\%$ ($3\sigma$) confidence.

\subsection{Differential algebra}\label{sec:da}

\Glsentrylong{da} is a computing technique based on the idea that it is possible to extract more information from a function $f$ than its mere value $y=f(x)$ at $x$. It replaces the algebra of \glsentrylong{fp} numbers with a new algebra of Taylor polynomials to compute the given order Taylor expansion of any function $f$ in $v$ variables that is $\mathcal{C}^{n+1}$ in the domain of interest $[-1,1]^v$~\cite{Berz1999}. The notation used is
\begin{equation}
    f \approx [f] = \mathcal{T}_f(\delta\vb*{x})
\end{equation}
where $\delta\vb*{x}=\{\delta x_1,\ldots,\delta x_v\}^T$ are the $v$ independent \glsentryshort{da} variables.

As for the \glsentrylong{fp} algebra, the four arithmetic operations ($+,-,\ast,\divisionsymbol$), elementary functions (e.g. trigonometric functions, powers, exponential and logarithm), derivation, integration, function composition and inversion are well defined in \glsentryshort{da}. The latter technique is particularly suitable for the solution of implicit equations and the computation of the Taylor expansion of the flow of differential equations in terms of their \glsentrylongpl{ic}.

\subsubsection{Polynomial bounding}

Consider the \glsentryshort{da} of order $n$ in $v$ variables, and a multivariate polynomial $p(\delta\vb*{x})$ with coefficients $a_k\in\mathbb{R}$. Define as \emph{odd} the coefficients for which at least one \glsentryshort{da} variable appears with an odd exponent (e.g. $\delta x_1$ or $\delta x_1^2\delta x_2$), and denote them with the superscript ``$o$". Likewise, define as \emph{even} coefficients those for which all variables appear with even exponents (e.g $\delta x_1^2$ or $\delta x_1^2\delta x_2^2$), and denote them with the superscript ``$e$". Since $\delta\vb*{x}\in[-1,1]^v$, an estimation of the range of variation, or bounds, of $p(\delta\vb*{x})$ can be obtained as
\begin{equation}
    \begin{aligned}
        b_l &= a_0 - \sum_{k=1}^{n_o}\abs{a^o_k} + \sum_{k=1}^{n_e} \min\{0, a^e_k\}\\
        b_u &= a_0 + \sum_{k=1}^{n_o}\abs{a^o_k} + \sum_{k=1}^{n_e} \max\{0, a^e_k\}
    \end{aligned}
    \label{eq:poly_bounds}
\end{equation}
with $b_l,b_u$ lower and upper bounds such that $p(\delta\vb*{x})\in[b_l,b_u]$, $a_0$ constant part and $n_o,n_e$ number of odd and even coefficients respectively. For first-order polynomials the bounds computed with \cref{eq:poly_bounds} are rigorous in the sense that they coincide with the extrema of $p(\delta\vb*{x})$ when evaluated in $[-1,1]^v$. In this case all even coefficients are identically zero and \cref{eq:poly_bounds} simplifies to
\begin{equation}
    p(\delta\vb*{x})\in \left[a_0\pm\sum_{k=1}^v\abs{a^o_k}\right]
\end{equation}

\subsubsection{Uncertainty propagation}

Given the nonlinear transformation $\vb*{f}:\mathbb{R}^n\to\mathbb{R}^m$ and two multivariate random variables $\vb*{X},\vb*{Y}$ with 
\glsentryfull{pdf} $p_{\vb*{X}},p_{\vb*{Y}}$ respectively, the \glsentryshort{up} problem requires the estimation of $p_{\vb*{Y}}$ from the only knowledge of $p_{\vb*{X}}$ and $\vb*{f}$. In \glsentryshort{da} framework, $\vb*{X}$ is represented as a vector of first-order multivariate polynomials in the form
\begin{equation}
    [\vb*{x}] = \bar{\vb*{x}} + \vb*{\beta}\odot\delta\vb*{x}
\end{equation}
with $\bar{\vb*{x}}=\mathbb{E}[\vb*{X}]$ the expected value of $\vb*{X}$, $\vb*{\beta}$ the semi-amplitude of the range of variation of $[\vb*{x}]$, $\delta\vb*{x}$ the first-order variations of $\vb*{X}$ around $\bar{\vb*{x}}$ and $\odot$ the Hadamard (or element-wise) product. The Taylor expansion of $\vb*{Y}$ in terms of $\delta\vb*{x}$ is then obtained evaluating $\vb*{f}$ in \glsentryshort{da} framework as
\begin{equation}
    [\vb*{y}] = \vb*{f}([\vb*{x}]) = \mathcal{T}_{\vb*{y}}(\delta\vb*{x})
    \label{eq:da_eval}
\end{equation}

In the context of \glsentryshort{up}, different methods were proposed to efficiently compute the statistical properties of $\vb*{Y}$ from $\mathcal{T}_{\vb*{y}}(\delta\vb*{x})$ thus avoiding the needs for expensive \glsentryfull{mc} simulations~\cite{Valli2013}.

\subsubsection{Low-order automatic domain splitting}\label{sec:loads}

Taylor polynomials are only a local approximation of the function $\vb*{f}$ around the expansion point, and a single polynomial might fail to accurately represent the uncertainty for $\abs{\delta\vb*{x}}\gg 0$. To overcome this limitation, \glsentryfull{ads} was firstly proposed by \cite{Wittig2015a} to monitor the accuracy of the expansion at the edges of its domain and recursively split the initial polynomial into smaller subdomains in which the required accuracy is satisfied.

Instead, a novel splitting technique specific for second-order Taylor polynomials is used in this work~\cite{Fossa2022,Losacco2023}. The algorithm, named \glsentryfull{loads}, defines a \glsentryfull{nli} computed from the Jacobian of the transformation $\partial\vb*{f}/\partial\vb*{x}$ that was inspired by the nonlinearity measure introduced by \cite{Junkins2004}. In this paper, $\vb*{f}$ is referred as the \emph{target} function of the \glsentryshort{loads} scheme.

Given an initial domain $[\vb*{x}^{(0)}]$, the algorithm starts by evaluating $\vb*{f}$ in $[\vb*{x}^{(0)}]$ and computing the \glsentryshort{nli} of the output $[\vb*{y}^{(0)}]$ as follows.

The Jacobian of the transformation is firstly computed as the sum of its constant part $\overline{\vb*{J}}$ and first-order expansion $\delta\vb*{J}(\delta\vb*{x})$ as
\begin{equation}
    [\vb*{J}] = \dfrac{\partial\vb*{f}}{\partial\vb*{x}} = \overline{\vb*{J}} + \delta\vb*{J}
\end{equation}
since $\delta\vb*{J}\equiv\vb*{0}$ for a linear transformation, a measure of nonlinearity can be obtained from the magnitude of $\delta\vb*{J}$. For each entry $\delta J_{ij}$, rigorous bounds are firstly obtained from \cref{eq:poly_bounds} as
\begin{equation}
    b_{ij}=\sum_{k=1}^n \abs{a_{ij,k}}
\end{equation}
where $a_{ij,k}$ are the first-order coefficients of $J_{ij}$, i.e. $\delta J_{ij}=\sum_{k=1}^n a_{ij,k}\delta x_k$ and thus $\delta J_{ij}\in[-b_{ij},b_{ij}]$. The \glsentrylong{nli} is then defined as
\begin{equation}
    \nu = \dfrac{\sqrt{\sum_{i=1}^m\sum_{j=1}^n\left(\sum_{k=1}^n \abs{a_{ij,k}}\right)^2}}{\sqrt{\sum_{i=1}^m\sum_{j=1}^n \overline{J}_{ij}^2}}
    \label{eq:nli_def}
\end{equation}
which is the ratio between the Frobenius norm of the matrix of the bounds and the Frobenius norm of the constant part of the Jacobian $\overline{\vb*{J}}$.

The computed \glsentryshort{nli} is then compared with an imposed error threshold $\varepsilon$ and a decision is made. If $\nu\leq\varepsilon$ the algorithm ends and $[\vb*{y}^{(0)}]$ is retained. On the contrary, if $\nu>\varepsilon$, the output is rejected and $[\vb*{x}^{(0)}]$ is split into three polynomials $[\vb*{x}^{(j)}]$ with each one covering $1/3$ of the initial domain. 

To perform the split, a second index is firstly computed to identify the direction that contributes the most to the nonlinearity in $[\vb*{y}]$. For each independent \glsentryshort{da} variable $\delta x_d,\ d\in[1,v]$, a directional Jacobian is built as the function composition of $[\vb*{J}]$ and $\delta\vb*{x}_d$ as
\begin{equation}
    [\vb*{J}_d] = [\vb*{J}]\circ\delta\vb*{x}_d
\end{equation}
where $\delta\vb*{x}_d$ is a vector of zeros except for the $d^{th}$ component set equal to $\delta x_d$, i.e. $\delta\vb*{x}_d=\{0,\ldots,\delta x_d,\ldots,0\}^T$, and $\circ$ denotes function composition. A directional \glsentryshort{nli} $\nu_d$ is then computed from $[\vb*{J}_d]$ as in \cref{eq:nli_def} and the splitting direction $d_{s}$ selected for the maximum $\nu_d$
\begin{equation}
    d_{s} = \argmax_d\{\nu_d\}
\end{equation}
Knowing $d_{s}$, the three subdomains $[\vb*{x}^{(j)}]$ for $j\in[1,3]$ are computed from $[\vb*{x}^{(0)}]$ as
\begin{equation}
    [\vb*{x}^{(j)}] = [\vb*{x}^{(0)}]\circ\left\{\delta x_1,\ldots,\frac{2}{3}\cdot(j-2)+\dfrac{1}{3}\cdot\delta x_{d_s},\ldots,\delta x_v\right\}^T
\end{equation}

The target function is then evaluated on the set $\{[\vb*{x}^{(j)}]\}$ and the procedure iterated until the accuracy is satisfied in each subdomain. The result is a set of polynomials for both the initial domain and its image through $\vb*{f}$. These sets are named \emph{manifolds} and denoted as $M_{\vb*{x}}=\{[\vb*{x}^{(i)}]\}$ and $M_{\vb*{y}}=\{[\vb*{y}^{(i)}]\}$ respectively where $i\in[1,N]$ and $N$ is the total number of domains after splitting.

\subsubsection{Domain merging}\label{sec:merging}

Consider a second function $\vb*{g}:\mathbb{R}^m\to\mathbb{R}^p$ evaluated on the manifold $M_{\vb*{y}}$ to obtain $M_{\vb*{z}}$. Nonlinearities in $[\vb*{z}^{(i)}]$ can be smaller than that of $[\vb*{y}^{(i)}]$ and one or more \glsxtrshortpl{nli} computed from $\partial\vb*{g}/\partial\vb*{y}$ can fall far below the splitting threshold, i.e. $\nu^{(i)}\ll\varepsilon$. In this case is thus conceivable to develop an algorithm to recombine the polynomials $[\vb*{z}^{(i)}]$ and reduce the size of $M_{\vb*{z}}$. The merging scheme used in this work stems from this observation and requires the knowledge of the entire \emph{splitting history} of each polynomial in $M_{\vb*{z}}$ to attempt a recombination of the different domains~\cite{Losacco2023}. This information comprises all directions $d_s$ and shifts $j$ along which $[\vb*{z}^{(i)}]$ was split and allows to identify domains that generated from the same parent and thus potentially mergeable.

Given the manifold $M_{\vb*{z}}$, the merging scheme proceeds as follows. Firstly, the domains are grouped by \emph{depth} of split, i.e. by the number of times they have been split. Secondly, starting from the sub-manifold with the greatest depth, the histories of its elements are compared to identify all triplets sharing a common parent. For each triplet, the domain of the parent is recovered as
\begin{equation}
    [\vb*{z}^{(i)}] = [\vb*{z}^{(2)}]\circ\left\{\delta x_1,\ldots,3\delta x_{d_s},\ldots,\delta x_v\right\}^T
\end{equation}
with $d_s$ last splitting direction and $[\vb*{z}^{(2)}]$ central polynomial. A decision on the potential merge is then taken after computing the \glsentryshort{nli} of $[\vb*{z}^{(i)}]$: if $\nu\leq\varepsilon$, $[\vb*{z}^{(i)}]$ is retained and added to the sub-manifold with a shallower depth. If $\nu>\varepsilon$, $[\vb*{z}^{(i)}]$ is discarded and the three domains of the original triplet added to the final manifold $M'_{\vb*{z}}$. This procedure is repeated until all sub-manifolds have been emptied. At each stage, if an incomplete triplet is found, no merge is attempted and its elements are directly added to $M'_{\vb*{z}}$.

\subsection{Orbit determination}\label{sec:od}
The problem of \glsentrylong{od} consists in estimating the state of an orbiting object given the measurements provided by one or more sensors. Among the different approaches to the \glsentryshort{od} problem, this paper exploits the so-called batch algorithms introduced in the followings.
\subsubsection{Batch orbit determination algorithms}\label{sec:batch_od}
Solving the batch \glsentryshort{od} problem requires to find an estimate of the orbit state $\hat{\vb*{x}}_0=\hat{\vb*{x}}(t_0)$, at some reference epoch $t_0$, that minimizes the residuals between the real observations $(\vb*{y}_k)_{k = 1}^{m}$ collected at epochs $(t_k)_{k = 1}^{m}$ and belonging to some observation space $\mathbf{Y} \subseteq \mathbb{R}^p$, and the expected (or theoretical) measurements $(\hat{\vb*{y}}_k)_{k = 1}^{m}$. The latter are computed by propagating $\hat{\vb*{x}}_0$ to $t_k$ and projecting the state $\hat{\vb*{x}}_k$ onto the observable space, i.e.,
\begin{equation}
    \hat{\vb*{y}}_k = \vb*{h}_k(\hat{\vb*{x}}_k) \qquad k\in[1,m]
\end{equation}
where $\vb*{h}_k$ is the observation function at epoch $t_k$. Assuming that the estimated state $\hat{\vb*{x}}_k$ is close to the true object's state $\vb*{x}_k$ at epoch $t_k$, $\vb*{h}_k$ can be linearized using a first-order Taylor expansion as
\begin{equation}
\begin{aligned}
    \vb*{h}_k(\vb*{x}_k) &\approx \vb*{h}_k(\hat{\vb*{x}}_k)+\tilde{\vb*{H}}_k\delta\vb*{x}_k\\
    &=\vb*{h}_k(\hat{\vb*{x}}_k)+\tilde{\vb*{H}}_k\vb*{\Phi}(t_k,t_0)\delta\vb*{x}_0\\
    &=\vb*{h}_k(\hat{\vb*{x}}_k)+\vb*{H}_k\delta\vb*{x}_0
\end{aligned}
\label{eq:obs_func}
\end{equation}
with $\vb*{\Phi}(t_k,t_0)$ the \glsentryfull{stm} from $t_0$ to $t_k$, $\tilde{\vb*{H}}_k$ the partials of $\vb*{h}_k$ with respect to $\vb*{x}_k$ and evaluated at $\hat{\vb*{x}}_k$, and $\vb*{H}_k$ the partials of $\vb*{h}_k$ with respect to $\vb*{x}_0$ and evaluated at $\hat{\vb*{x}}_0$. The last two matrices are defined as
\begin{subequations}
    \begin{gather}
        \tilde{\vb*{H}}_k = \left.\dfrac{\partial\vb*{h}_k}{\partial\vb*{x}_k}\right|_{\hat{\vb*{x}}_k}\\
        \vb*{H}_k = \tilde{\vb*{H}}_k\vb*{\Phi}(t_k,t_0) = \left.\dfrac{\partial\vb*{h}_k}{\partial\vb*{x}_k}\dfrac{\partial\vb*{x}_k}{\partial\vb*{x}_0}\right|_{\hat{\vb*{x}}_0}\label{eq:obs_jac}
    \end{gather}
\end{subequations}
The residuals are then readily computed as
\begin{equation}
\begin{aligned}
    \vb*{e} &= \vb*{y} - \hat{\vb*{y}}\\
    &= \vb*{y} - \vb*{h}(\hat{\vb*{x}}) - \vb*{H}\delta\vb*{x}_0\\
    &= \Delta\vb*{y} - \vb*{H}\delta\vb*{x}_0
\end{aligned}
\end{equation}
where the total observations vectors and Jacobian are defined as
\begin{subequations}
\begin{align}
    \vb*{y} &= [\vb*{y}^T_1,\ldots,\vb*{y}^T_m]^T\\
    \hat{\vb*{y}} &= [\hat{\vb*{y}}^T_1,\ldots,\hat{\vb*{y}}^T_m]^T\\
    \vb*{h}(\hat{\vb*{x}}) &= [\vb*{h}_1(\hat{\vb*{x}}_1),\ldots,\vb*{h}_k(\hat{\vb*{x}}_k)]^T\\
    \vb*{H} &= [\vb*{H}_1^T,\ldots,\vb*{H}_m^T]^T
\end{align}
\end{subequations}
The three vectors $\vb*{y},\hat{\vb*{y}}$ and $\vb*{h}(\hat{\vb*{x}})$ have size $mp\times 1$ while the matrix $\vb*{H}$ has size $mp\times n$, where $n$ is the dimension of $\hat{\vb*{x}}_k$.

Given an initial guess for $\hat{\vb*{x}}_0$ denoted as $\hat{\vb*{x}}_0^0$, the objective of the \glsentryshort{od} algorithm is to iteratively solve for $\delta\vb*{x}_0$ and compute an updated estimate $\hat{\vb*{x}}_0^j$ as
\begin{equation}
    \hat{\vb*{x}}_0^j = \hat{\vb*{x}}_0^{j-1} + \delta\vb*{x}_0
\end{equation}
with $j$ iteration number. The recursive procedure terminates as soon as one of the following criteria is satisfied
\begin{itemize}
    \setlength{\itemsep}{0pt plus 1pt}
    \item $\norm{\hat{\vb*{y}}^j-\hat{\vb*{y}}^{j-1}}_2<\varepsilon_{res}$, i.e. the residuals at two subsequent iterations differ by less than the residual tolerance $\varepsilon_{res}$.
    \item $\abs{J^j-J^{j-1}}<\varepsilon_{opt}$, i.e. the cost functions at two subsequent iterations differ by less than the optimality tolerance $\varepsilon_{opt}$.\footnote{the cost function $J$ is defined in \cref{sec:ls,sec:lsar}.}
    \item $\norm{\delta\hat{\vb*{x}}_0}_2<\varepsilon_{step}$, i.e. the state update is smaller than the step tolerance $\varepsilon_{step}$.
    \item $j=j_{max}$, i.e. a maximum number of iterations has been reached.
\end{itemize}

\subsubsection{Least squares}\label{sec:ls}
The \glsentryshort{ls} solution to the \glsentryshort{od} problem solves for $\delta\vb*{x}_0$ that minimizes the residuals in the \glsentryshort{ls} sense. The cost function is given by
\begin{equation}
    J_{LS} = \frac{1}{2}\vb*{e}^T\vb*{W}\vb*{e} = \frac{1}{2}\left[\Delta\vb*{y}-\vb*{H}\delta\vb*{x}_0\right]^T\vb*{W}\left[\Delta\vb*{y}-\vb*{H}\delta\vb*{x}_0\right]
    \label{eq:cost_fun_ls}
\end{equation}
with $\vb*{W} = \text{diag}(1/\sigma_i^2)$ the diagonal weight matrix where $\sigma_i$ is the sensor noise standard deviation corresponding to the $i^{th}$ observation. The state update that minimizes \cref{eq:cost_fun_ls} is obtained using the \glsentrylong{lm} algorithm~\cite{More1980} (damped \glsentryshort{ls}) as
\begin{equation}
    \delta\vb*{x}_0 = (\vb*{H}^T\vb*{W}\vb*{H} + \lambda \vb*{I})^{-1}\vb*{H}^T\vb*{W}\Delta\vb*{y}
    \label{eq:lm_iter}
\end{equation}
with $\lambda$ the damping parameter and $I$ the identity matrix. The covariance matrix of $\hat{\vb*{x}}_0$ is finally estimated as
\begin{equation}
    \vb*{P}_0 = (\vb*{H}^T\vb*{W}\vb*{H})^{-1}
\end{equation}

\subsubsection{Minimum \texorpdfstring{$L_1$}{L1} norm estimates}\label{sec:lsar}
The minimum $L_1$ norm solution to the \glsentryshort{od} problem was proven to be a robust alternative to \glsentryshort{ls} optimization capable of rejecting outlier measurements in the presence of redundant observations~\cite{Narula1982,Portnoy1997,Branham2005,Gehly2016,Prabhu2021}. In this work, the \glsentryfull{lsar} algorithm proposed in~\cite{Prabhu2021} is used to solve the \glsentryshort{od} problem in a $L_1$ sense. The cost function to be minimized is
\begin{equation}
    J_{LSAR} = \sum_{i=1}^{mp}\tilde{w}_i\abs{e_i}
    \label{eq:cost_fun_lsar}
\end{equation}
with $e_i$ the components of $\vb*{e}$ and $\tilde{w}_i$ the weights defined as $\tilde{w}_i=1/(1.24\sigma_i)$~\cite{Rosenberg1977}. The problem is recast as a \glsentryfull{lpp} in the form
\begin{equation}
    \min_{\vb*{z}} \vb*{c}^T\vb*{z} \qquad \text{subject to}\ \vb*{Q}\vb*{z}\leq\vb*{k}
    \label{eq:lsar_lpp}
\end{equation}
with
\begin{subequations}
\begin{align}
    \vb*{Q} &= \begin{bmatrix} -\vb*{H} & -\vb*{I}_{mp\times mp} \\ \vb*{H} & -\vb*{I}_{mp\times mp} \end{bmatrix}\\
    \vb*{z} &= \begin{bmatrix} \delta\vb*{x}_0^T & \vb*{s}^T \end{bmatrix}^T\\
    \vb*{k} &= \begin{bmatrix} -\Delta\vb*{y}^T & \Delta\vb*{y}^T \end{bmatrix}^T\\
    \vb*{c} &= \begin{bmatrix} \vb*{0}_{6\times 1}^T & \tilde{\vb*{w}}^T \end{bmatrix}^T
\end{align}
\end{subequations}
and $\vb*{s}$ slack variables for the constraint vector $\vb*{e}$ such that
\begin{equation}
    \abs{e_i}\leq s_i \qquad i=1,\ldots,mp
\end{equation}

Throughout this work, the solver MOSEK\footnote{\href{https://www.mosek.com}{https://www.mosek.com}} with its interior point algorithm is used to efficiently solve \cref{eq:lsar_lpp}. The covariance matrix of $\hat{\vb*{x}}_0$ is then estimated as
\begin{equation}
    \tilde{P}_0 = (\vb*{H}^T\tilde{\vb*{W}}\vb*{H})^{-1}
\end{equation}
with $\tilde{\vb*{W}}=\text{diag}(\tilde{w}_i^2)=\vb*{W}/1.24^2$.

\section{Robust orbit determination}\label{sec:robust_od}

The full batch \glsentrylong{od} chain is presented in this section. Given the \glsentryfull{so} of interest, an estimate of its orbit state at epoch $t_0$ is firstly obtained with the \glsentryshort{da}-based \glsentryfull{iod} routine described in \cref{sec:iod}. This estimate is then propagated forward in time using \glsentrylong{mf} techniques to minimize the computational effort and its uncertainty progressively reduced with a custom pruning scheme as soon as new measurements become available. Potential outliers in the input measurements are also identified and discarded in this step. A batch \glsentryshort{od} algorithm is finally run on the pruned state and retained observations to compute the final estimate. Both \glsentryshort{ls} and \glsentryshort{lsar} techniques are tested and their performance discussed in \cref{sec:application}. The overall workflow is summarized in \cref{fig:flowchart}. This paper assumes observations from optical telescopes, but each tool can be easily adapted to other types of measurements. For the \glsentryshort{iod} algorithm in \cref{sec:iod} it is already the case with a recently submitted paper that extends the algorithm to radar measurements~\cite{Fossa2023}.

\subsection{Initial orbit determination}\label{sec:iod}

Consider a set of $N$ tuples of angular measurements as provided by a ground based optical sensor while observing an unknown \glsentryshort{so}
\begin{equation}
    \left\{t_i;\left(\alpha_i;\sigma_{\alpha_i}\right),\left(\delta_i;\sigma_{\delta_i}\right)\right\}\qquad i\in[1,N]
    \label{eq:meas_topo_set}
\end{equation}
with $\alpha_i,\delta_i$ topocentric right ascension and declination of the \glsentryshort{so} at epoch $t_i$ and $\sigma_{\alpha_i},\sigma_{\delta_i}$ associated standard deviations of the sensor noise, assumed as uncorrelated Gaussian white noise as described in \cref{sec:meas_model}.

The goal of the \glsentryshort{iod} is to estimate the state of the object from the available set of measurements. The problem is well-known in literature and typically solved by processing the nominal angular measurements only~\cite{Laplace1780,Gauss1809}. Instead, the approach here exploited makes use of \glsentryshort{da} to obtain an analytic map between the \glsentrylongpl{ci} of the processed measurements and the uncertainty of the computed solution. It was firstly developed under the assumption of Keplerian dynamics by~\cite{Pirovano2020e} for optical measurements and~\cite{Losacco2023a} for radar sensors, and recently extended in~\cite{Fossa2023} to accommodate any dynamical model for the target object. The method sets up an iterative procedure to estimate the topocentric ranges of the object at the epoch of the first, middle, and last available measurements, namely $\rho_k=\{\rho_1,\rho_2,\rho_3\}$. The initial guess for these ranges is provided by the classical Gauss solution~\cite{Gauss1809}. Two Lambert's problems~\cite{Izzo2015} are then solved between $t_1,t_2$ and $t_2,t_3$ respectively. A nonlinear solver is then used to update $\rho_j$ such that the velocity vectors at the intermediate epoch coincide. An estimate of the object state is then directly available from this solution.

A polynomial expansion of the state with respect to the observables can be obtained by initializing the processed measurements as Taylor variables and performing all operations in the \glsentryshort{da} framework. Under the assumption of uncorrelated measurements corrupted by Gaussian white noise, these are written as
\begin{equation}
\begin{aligned}
    [\vb*{\alpha}] &= \vb*{\alpha}+c\vb*{\sigma}_{\vb*{\alpha}}\odot\delta\vb*{\alpha}\\
    [\vb*{\delta}] &= \vb*{\delta}+c\vb*{\sigma}_{\vb*{\delta}}\odot\delta\vb*{\delta}
\end{aligned}
\label{eq:iod_alpha_delta_da}
\end{equation}
where $c>0$ is the $z$-score introduced in \cref{sec:meas_model} and
\begin{equation}
\begin{matrix}
    \vb*{\alpha} = [\alpha_1\ \alpha_2\ \alpha_3]^T & \vb*{\sigma}_{\vb*{\alpha}} = [\sigma_{\alpha_1}\ \sigma_{\alpha_2}\ \sigma_{\alpha_3}]^T \\
    \vb*{\delta} = [\delta_1\ \delta_2\ \delta_3]^T & \vb*{\sigma}_{\vb*{\delta}} = [\sigma_{\delta_1}\ \sigma_{\delta_2}\ \sigma_{\delta_3}]^T
\end{matrix}
\end{equation}
The state vectors at $t_j$ can then be written as
\begin{equation}
    [\vb*{x}_j]=\mathcal{T}_{\vb*{x}_j}\left(\delta\vb*{\alpha},\delta\vb*{\delta}\right) \qquad j\in[1,3]
    \label{eq:iod_sol_kep}
\end{equation}

The solution in \cref{eq:iod_sol_kep} is obtained under the assumption of unperturbed Keplerian motion (since based on the recursive solution of two Lambert problems) and might be inaccurate for long observation windows in which deviations from the nominal two-body trajectory cannot be neglected. A \glsentryshort{da}-based shooting scheme is thus employed to account for the effects of the $J_2$ term of the Earth's gravitational potential. Given the needs for a fast propagation method at this stage, the analytical formulation of the $J_2$-perturbed dynamics proposed in \cite{Armellin2018h} is employed. Starting from \cref{eq:iod_sol_kep}, an updated estimate is then obtained as
\begin{equation}
    [\vb*{x}_j^{J_2}]=\mathcal{T}_{\vb*{x}_j^{J_2}}\left(\delta\vb*{\alpha},\delta\vb*{\delta}\right) \qquad j\in[1,3]
    \label{eq:iod_sol_j2}
\end{equation}

To maintain an accurate polynomial representation of the state in the presence of large measurements errors, the \glsentryshort{iod} algorithm is wrapped within the \glsentryshort{loads} framework introduced in \cref{sec:loads}. A manifold of polynomials $M_{\vb*{x}_j^{J_2}}$ is thus automatically generated while computing the expansion of the \glsentryshort{iod} solution with respect to $(\vb*{\alpha},\vb*{\delta})$ such that the nonlinear transformation is accurately described across its entire domain. These manifolds are defined as
\begin{equation}
    M_{\vb*{x}_j^{J_2}} = \left\{[\vb*{x}_j^{J_2}]^{(i)}\right\} \qquad i\in[1,N],\ j\in[1,3]
    \label{eq:iod_sol_j2_loads}
\end{equation}
where $N$ is the number of domains generated by the \glsentryshort{loads} algorithm and $j$ denotes the observation epoch.

\subsection{Multifidelity uncertainty propagation}\label{sec:mf_up}

\Glsentrylong{mf} methods are introduced to improve the computational efficiency of the prediction step while guaranteeing the required accuracy for the final estimate. In particular, a bifidelity technique based on two distinct dynamical models is employed.

The \glsentrylong{lf} model is the \glsentryfull{sgp} model used by \glsentryshort{norad} to produce daily \glsentryfullpl{tle} of all tracked \glsentryshortpl{so}~\cite{Vallado2006a}. Given its analytical formulation, \glsentryshort{sgp} is easily embedded in the \glsentryshort{da}-\glsentryshort{loads} framework to efficiently propagate orbit states and associated uncertainties as Taylor polynomials. A \glsentrylong{hf} correction is yet needed to retain the accuracy of the final solution. A numerical propagator is thus used to integrate the \glsentrylong{hf} orbital dynamics as well as the effects of \glsentrylong{pn} via \cref{eqn:ode_snc}. In this work, $\vb*{f}(\vb*{x},t)$ are the Gauss variational equations which describe the motion of \glsentryshortpl{so} in the perturbed Keplerian dynamics~\cite{Battin1999}. A representation of the orbit state at epoch $t$ is then maintained through a manifold of Taylor polynomials $M_{MF}(t)=\{[\vb*{x}_{MF}^{(i)}(t)]\}$ and a manifold of covariance matrices $M_{PN}(t)=\{\vb*{P}_{PN}^{(i)}(t)\}$ characterized by a one-to-one correspondence between the elements in $M_{MF}(t)$ and the elements in $M_{PN}(t)$. The polynomials $[\vb*{x}_{MF}^{(i)}(t)]$ represent the initial uncertainty on the \glsentryshort{iod} solution propagated at time $t$ under a deterministic dynamical model and are centered on the \glsentrylong{hf} reference trajectories solution of the first of \cref{eqn:ode_snc}. The covariance matrices $\vb*{P}_{PN}^{(i)}(t)$ represent instead the contribution of process noise effects and are obtained by numerically integrating the second of \cref{eqn:ode_snc}. At initial time $t_0$ the first manifold $M_{MF}(t_0)$ is thus given by \cref{eq:iod_sol_j2} for $j=1$ while all covariance matrices are initialized to zero, i.e. $\vb*{P}_{PN}^{(i)}(t_0) = \vb*{0}\ \forall i$. Note that, at any time $t$, the polynomials $[\vb*{x}_{MF}^{(i)}(t)]$ are function of the six independent \glsentryshort{da} variables $(\delta\vb*{\alpha},\delta\vb*{\delta})$, i.e.
\begin{equation}
    [\vb*{x}_{MF}^{(i)}(t)]=\mathcal{T}_{\vb*{x}_{MF}^{(i)}(t)}(\delta\vb*{\alpha},\delta\vb*{\delta}) \qquad \forall i
\end{equation}

Consider the propagation time span $[t_{k-1},t_k]$ between two subsequent observations, and denote with $M_{MF}(t_{k-1})$ and $M_{PN}(t_{k-1})$ the manifolds of Taylor polynomials and covariance matrices at time $t_{k-1}$ respectively. The solution manifolds $M_{MF}(t_{k}), M_{PN}(t_{k})$ at time $t_k$ are then computed as follows. Firstly, $M_{MF}(t_{k-1})$ is processed with the \glsentryshort{loads} scheme setting \glsentryshort{sgp} as target function $\vb*{f}$. The results is a \glsentrylong{lf} solution manifold $M_{LF}(t_{k})$ whose size is generally different from that of $M_{MF}(t_{k-1})$. In the second step, the centers of the polynomials $[\vb*{x}^{(i)}(t_{k-1})]$ that maps to $[\vb*{x}^{(i)}(t_{k})]$ in $M_{LF}(t_{k})$ are propagated \glsentrylong{pw} in \glsentrylong{hf} using \cref{eqn:ode_snc} to obtain the reference solutions $\{\vb*{x}^{(i)}_{HF}(t_{k})\}$ and covariance matrices $\{\vb*{P}^{(i)}_{PN}(t_{k})\}$ for each polynomial in $M_{LF}(t_{k})$. Note that \cref{eqn:ode_snc} is now subject to \glsentrylongpl{ic} $\bar{\vb*{x}}^{(i)}_{MF}(t_{k-1}),\vb*{P}^{(i)}_{PN}(t_{k-1})$, while $\vb*{P}^{(i)}_{PN}(t_{k})$ is the result of the effects of \glsentrylong{pn} and do not include the contribution of the initial uncertainty which was already taken into account in the \glsentrylong{lf} step. The \glsentrylong{mf} solution is then computed re-centering the Taylor expansions in $M_{LF}(t_{k})$ on the \glsentrylong{hf} trajectories $\{\vb*{x}^{(i)}_{HF}(t_{k})\}$. The $i^{th}$ domain of $M_{MF}(t_k)$ is thus given by
\begin{equation}
    [\vb*{x}_{MF}^{(i)}(t_{k})] = \vb*{x}^{(i)}_{HF}(t_{k}) + \left\{[\vb*{x}_{LF}^{(i)}(t_{k})] - \bar{\vb*{x}}_{LF}^{(i)}(t_{k})\right\}
    \label{eq:mf_state_no_noise}
\end{equation}
which is a Taylor polynomial centered at $\vb*{x}^{(i)}_{HF}(t_{k})$ and whose non-constant coefficients model the uncertainty due to the \glsentryshortpl{ic} only. Information on the cumulative effects of stochastic accelerations is kept separately in $M_{PN}(t_k)$, and added to \cref{eq:mf_state_no_noise} when computing the polynomial bounds in observables space as described in \cref{sec:proj}. These two manifolds are collectively referred as $M_{\vb*{x}}(t_k)=\{([\vb*{x}_{MF}^{(i)}(t_{k})],\vb*{P}_{PN}^{(i)}(t_{k}))\}$ and constitute all the information on the \glsentryshort{so} state that is retained between subsequent steps.

For computational purposes, alternate equinoctial elements are used for propagation. This set of coordinates is defined in terms of the six Keplerian parameters as~\cite{Horwood2011}
\begin{equation}
    \begin{aligned}
        n &= \sqrt{\mu/a^3}\\
        f &= e\cos(\omega + \Omega)\\
        g &= e\sin(\omega + \Omega)\\
        h &= \tan(i/2)\cos(\Omega)\\
        k &= \tan(i/2)\sin(\Omega)\\
        \lambda &= \Omega + \omega + M
    \end{aligned}
    \label{eq:alt_el}
\end{equation}
with $a$ the semi-major axis, $e$ the eccentricity, $i$ the inclination, $\Omega$ the \glsentrylong{raan}, $\omega$ the \glsentrylong{aop}, and $M$ the mean anomaly, where $n$ is the mean motion, and $\lambda$ the mean longitude. Given the quasi-linearity of the \glsentryshortpl{ode} governing their evolution over time in perturbed Keplerian dynamics, this choice leads to smaller number of domains generated in the \glsentrylong{lf} step compared to either Cartesian, Keplerian or (modified) equinoctial parameters~\cite{Fossa2022c,Fossa2022b}.

\subsection{Projection onto observables space}\label{sec:proj}

Consider a set of osculating orbital elements $\vb*{x}(t_k)$ at epoch $t_k$ and an optical telescope as described in \cref{sec:meas_model}. The spherical coordinates $(\rho_k,\alpha_k,\delta_k)$ of the \glsentryshort{so} as seen from the telescope at epoch $t_k$ are given by \cref{eq:cart2sph} with $\vb*{\rho}=[\rho_x\ \rho_y\ \rho_z]^T$ \glsentryshort{los} vector obtained from \cref{eq:so_pos_from_range}. Moreover, the \glsentryshort{so}'s position $\vb*{r}_k$ is easily obtained by transforming $\vb*{x}(t_k)$ to Cartesian parameters. These operations are carried out in the \glsentryshort{da} framework using the \glsentryshort{loads} algorithm to accurately map the uncertainties in $\vb*{x}(t_k)$ through the sequence of nonlinear transformations. However, the manifold $M_{MF}(t_k)=\{[\vb*{x}_{MF}^{(i)}(t_k)]\}$ output of \cref{sec:mf_up} cannot be directly processed if the effects of process noise on the projected states have to be taken into account. The last are in fact represented by the covariance matrices in $M_{PN}(t_k)$, and an intermediate manifold must be built to include this information. Note that each polynomial in $M_{MF}(t_k)$ is function of the independent \glsentryshort{da} variables $(\delta\vb*{\alpha},\delta\vb*{\delta})$ while there is no concept of Taylor expansion in $\vb*{P}_{PN}^{(i)}(t_k)$. Six additional independent variables $\delta\vb*{x}_k$ are thus considered to describe $\vb*{P}_{PN}^{(i)}(t_k)$ in the \glsentryshort{da} framework as
\begin{equation}
    [\vb*{x}_{PN}^{(i)}(t_k)] = \vb*{V}^{(i)}\cdot\left(c\sqrt{\vb*{\lambda}^{(i)}}\delta\vb*{x}_k\right)
    \label{eq:pn_cov_da}
\end{equation}
with $c$ the $z$-score introduced in \cref{sec:meas_model}, $\vb*{P}^{(i)}_{PN}(t_k)=\vb*{V}^{(i)}\vb*{\Lambda}^{(i)}\left(\vb*{V}^{(i)}\right)^T$ the covariance matrix eigendecomposition, and $\vb*{\lambda}^{(i)}$ the main diagonal of $\vb*{\Lambda}^{(i)}$. For each element in $M_{\vb*{x}}(t_k)=\{([\vb*{x}_{MF}^{(i)}(t_k)],\vb*{P}_{PN}^{(i)}(t_k))\}$ a Taylor expansion that represents the uncertainty due to both the \glsentryshortpl{ic} and the stochastic accelerations is thus given by
\begin{equation}
    \begin{aligned}
            [\vb*{x}_{MF,PN}^{(i)}(t_k)] &= [\vb*{x}_{MF}^{(i)}(t_k)] + [\vb*{x}_{PN}^{(i)}(t_k)]\\
            &= \mathcal{T}_{\vb*{x}_{MF,PN}^{(i)}(t_k)}(\delta\vb*{\alpha},\delta\vb*{\delta},\delta\vb*{x}_k)
    \end{aligned}
    \label{eq:poly_mf_pn}
\end{equation}
which is a polynomial in twelve independent \glsentryshort{da} variables. These polynomials are then collected into a manifold $M_{MF,PN}(t_k)$ which is subsequently projected onto the observables space $(\alpha,\delta)$. The result is a manifold of expected measurements $M_{\hat{\vb*{y}}}(t_k)$ with entries
\begin{equation}
    [\hat{\vb*{y}}^{(i)}(t_k)] = 
    \begin{bmatrix}
        [\hat{\alpha}^{(i)}(t_k)]\\
        [\hat{\delta}^{(i)}(t_k)]\\    
    \end{bmatrix}
    = \mathcal{T}_{\hat{\vb*{y}}^{(i)}(t_k)}(\delta\vb*{\alpha},\delta\vb*{\delta},\delta\vb*{x}_k)
    \label{eq:poly_proj}
\end{equation}
note that the cardinality of $M_{\hat{\vb*{y}}}(t_k)$ is generally greater that that of $M_{MF,PN}(t_k)$ since several splits might occur while evaluating the transformations in the \glsentryshort{da} framework. The domains $[\vb*{x}_{MF,PN}^{(i)}(t_k)]$ that maps onto $[\hat{\vb*{y}}^{(i)}(t_k)]$ are then easily reconstructed a posteriori to obtain a bijection between the elements of the two manifolds. This step is essential to correctly perform the subsequent domain pruning.

\subsection{Domain pruning}\label{sec:pruning}

Consider the manifold of observables $M_{\hat{\vb*{y}}}(t_k)$ estimated from $M_{MF,PN}(t_k)$ as in \cref{sec:proj} and the real measurements $(\alpha_k,\delta_k)$ taken from the same telescope at epoch $t_k$. These measurements are assumed to be corrupted by uncorrelated white noise with standard deviations $\sigma_{\alpha_k},\sigma_{\delta_k}$ respectively. The pruning algorithm seeks to reduce the size of $M_{\hat{\vb*{y}}}(t_k)$ exploiting the new information obtained from the sensor. The real measurements are firstly initialized as a \glsentryshort{da} vector
\begin{equation}
    [\vb*{y}] = \begin{bmatrix}\alpha_k\\ \delta_k\end{bmatrix} + c\cdot\begin{bmatrix}\sigma_{\alpha_k}\\ \sigma_{\delta_k}\end{bmatrix}\odot\delta\vb*{y}
    \label{eq:real_meas_da}
\end{equation}
with $c$ the $z$-score as in \cref{eq:ang_meas_da} and $\delta\vb*{y}$ the first-order deviations in $\alpha_k,\delta_k$. Rigorous bounds on $[\vb*{y}]$ are then given by the square interval $[\alpha_k\pm c\cdot\sigma_{\alpha_k}]\times[\delta_k\pm c\cdot\sigma_{\delta_k}]$ which is a simplification of \cref{eq:poly_bounds} for first-order univariate polynomials.

For each domain in $M_{\hat{\vb*{y}}}(t_k)$, its bounds are then estimated via \cref{eq:poly_bounds} and compared with the ones above. If the two intervals overlap, the domain is retained for the subsequent steps. The polynomial is otherwise discarded from $M_{\hat{\vb*{y}}}(t_k)$ since it is unlikely to include the true, yet unknown, orbit state. Similarly, if no intersection is found between $M_{\hat{\vb*{y}}}(t_k)$ and $[\vb*{y}]$, the measurement is deemed to be an outlier and no pruning is performed. The observation is subsequently excluded from the set processed with the batch \glsentryshort{od} algorithm described in \cref{sec:batch_od}. As last step, pruning is also performed on $M_{MF}(t_k)$ and $M_{PN}(t_k)$ so as to reestablish a bijection between the elements of the three sets before processing the next measurement.

\subsection{Combined prediction and pruning}\label{sec:seq_algo}

The overall algorithm for the prediction and pruning of orbit states is presented in this section. Starting from the \glsentryshortpl{ic} given by \cref{eq:iod_sol_j2_loads}, a sequential procedure is implemented to progressively tighten the bounds on $(\delta\vb*{\alpha},\delta\vb*{\delta})$ by exploiting the information available from subsequent measurements.

The estimated Cartesian state is firstly transformed into alternate equinoctial elements. This transformation is evaluated in the \glsentryshort{da} framework using the \glsentryshort{loads} scheme and the transformed \glsentryshortpl{ic} are thus represented as a manifold of polynomials $M_{MF}(t_0)$ function of $(\delta\vb*{\alpha},\delta\vb*{\delta})$. The state covariances due to process noise are also initialized as $M_{PN}(t_0)=\{\vb*{0}_{6\times 6}\}$ for each polynomial in $M_{MF}(t_0)$, and the two manifolds grouped into $M_{\vb*{x}}(t_0)=\{([\vb*{x}^{(i)}(t_0)],\vb*{0}_{6\times 6}^{(i)})\}$.

The measurement collected at time $t_k$ is denoted by $(\alpha_k,\delta_k)$, $1 \leq k \leq M$, and the associated $1\sigma$ uncertainties by $(\sigma_{\alpha_k},\sigma_{\delta_k})$, such that the developed algorithm requires $M$ sequential steps to estimate $M_{\vb*{x}}(t_{M})$ from $M_{\vb*{x}}(t_{0})$ while exploiting the information available from the measurements. Each step $k$ is further subdivided into the following operations

\begin{enumerate}
    \item Propagation of $M_{\vb*{x}}(t_{k-1})$ to the next epoch $t_k$ using the \glsentrylong{mf} scheme in \cref{sec:mf_up} to obtain $M_{\vb*{x}}(t_{k})$
    \item Construction of the inflated manifold $M_{MF,PN}(t_{k})$ and projection of the last onto the observables space $(\rho,\alpha,\delta)$ as in \cref{sec:proj} to obtain $M_{\hat{\vb*{y}}}(t_{k})$
    \item Pruning of $M_{\hat{\vb*{y}}}(t_{k})$ and $M_{\vb*{x}}(t_{k})$ with information from $(\alpha_k,\delta_k)$ as in \cref{sec:pruning} to obtain the reduced manifolds $M'_{\hat{\vb*{y}}}(t_{k}), M'_{\vb*{x}}(t_{k})$ as well as isolate potential outliers. If an outlier is found, no pruning is performed and $M'_{\vb*{x}}(t_{k})$ coincides with $M_{\vb*{x}}(t_{k})$
    \item Merging of the domains in $M'_{\vb*{x}}(t_{k})$ as in \cref{sec:merging} since nonlinearities are weaker in orbital elements space than in $(\rho,\alpha,\delta)$ and fewer domains are thus required to accurately capture the uncertainty in this space
\end{enumerate}

The procedure thus ends at epoch $t_{M}$ with the processing of the last available measurement. Two lists of indexes $K_{c}$ and $K_{o}$ corresponding to correlated and outlier measurements are also available at this stage.

\subsection{Differential algebra orbit determination}\label{sec:da_od}
Outputs of \cref{sec:seq_algo} include a manifold of osculating elements $M'_{\vb*{x}}(t_M)$ at the last measurement epoch $t_M$ and a set of correlated measurements $(\vb*{y}_k)_{k\in K_c}$. From these data a batch \glsentryshort{od} problem can be then set up as follows.

Firstly, the domains $M'_{\vb*{x}}(t_0)$ that map onto the pruned manifold at $t_M$ are reconstructed by applying the \emph{splitting history} of $M'_{\vb*{x}}(t_M)$ to the \glsentryshortpl{ic} $M_{\vb*{x}}(t_0)$. Among all domains in $M'_{\vb*{x}}(t_0)$, the center of the one characterized by the smallest residuals with respect to the correlated observations $(\vb*{y}_k)_{k\in K_c}$ is used as initial guess for the subsequent batch \glsentryshort{od}. The last is denoted as $\hat{\vb*{x}}_0^0$ where the superscript $``0"$ refers to the zeroth iteration of the algorithm. The best estimate $\hat{\vb*{x}}_0$ and associated covariance $\vb*{P}_0$ are then computed using either the \glsentryshort{ls} or \glsentryshort{lsar} algorithms described in \cref{sec:ls,sec:lsar} respectively. Within each iteration $j$ of the \glsentryshort{od} scheme, \glsentryshort{da} is again leveraged to automate the computation of the observation function $\vb*{h}_k$ defined in \cref{eq:obs_func} and its partials $\vb*{H}_k$ with respect to the current estimate $\hat{\vb*{x}}_0^j$. To do so, $\hat{\vb*{x}}_0^j$ is firstly initialized as the \glsentryshort{da} vector $[\hat{\vb*{x}}_0^j]=\hat{\vb*{x}}_0^j+\delta\vb*{x}_0$. The propagation to $t_k$ and projection onto the observables space are then carried out in the \glsentryshort{da} framework after setting the expansion order equal to one. The Taylor expansion of the theoretical measurements $\hat{\vb*{y}}_k$ at epoch $t_k$ is thus available as
\begin{equation}
    [\hat{\vb*{y}}_k] = \mathcal{T}_{\hat{\vb*{y}}_k}(\delta\vb*{x}_0)
    \label{eq:exp_y}
\end{equation}
$\vb*{h}_k(\hat{\vb*{x}}_k)$ and $\vb*{H}_k$ are then retrieved as the constant part and the linear part of \cref{eq:exp_y}, respectively. This avoids the needs for an analytical expression of $\partial\vb*{h}_k/\partial\vb*{x}_k$ and the integration of the $6\times 6$ variational equations for $\vb*{\Phi}$.

\begin{figure}[H]
    \centering
    \vspace{-2cm}
    \begin{tikzpicture}[node distance=2cm]
    
    \node (start) [startstop] {Start};
    \node (input) [io, below of=start, yshift=0.5cm] {Measurements $(\vb*{y}_k)_{k=1}^M$};
    \node (iod) [process, below of=input] {Perform \glsentryshort{iod} to obtain $M_{\vb*{x}}(t_k)$ for $k=1$};
    \node (proj) [process, below of=iod] {Project $M_{\vb*{x}}(t_k)$ to $M_{\hat{\vb*{y}}}(t_k)$};
    \node (prune) [process, below of=proj] {Prune $M_{\hat{\vb*{y}}}(t_k)$ and $M_{\vb*{x}}(t_k)$ to obtain $M'_{\vb*{x}}(t_k)$};
    \node (outlier) [decision, below of=prune, yshift=-1cm] {$M'_{\vb*{x}}(t_k)=\varnothing$};
    \node (merge) [smallprocess, below of=outlier, yshift=-0.75cm] {Merge $M'_{\vb*{x}}(t_k)$};
    \node (restore) [smallprocess, left of=merge, xshift=-2cm] {Restore $M_{\vb*{x}}(t_k)$ and discard $\vb*{y}_k$};
    \node (loop) [decision, below of=merge] {$k<M$};
    \node (prop) [process, below of=loop, yshift=-0.5cm] {Propagate $M_{\vb*{x}}(t_{k})$ to $t_{k+1}$};
    \node (count) [smallprocess, left of=prune, xshift=-4cm] {Set $k=k+1$};
    
    \node (repl) [process, right of=outlier, xshift=4cm] {Split $M_{\vb*{x}}(t_1)$ with pattern of $M'_{\vb*{x}}(t_M)$};
    \node (extract) [process, below of=repl] {Extract $\hat{\vb*{x}}_0^0$ and\\ retained $(\vb*{y}_k)_{k=1}^m$};
    \node (od) [process, below of=extract] {Perform \glsentryshort{od} on $\hat{\vb*{x}}_0^0$ and $(\vb*{y}_k)_{k=1}^m$};
    \node (out) [io, below of=od] {\Glsentryshort{od} solution $\hat{\vb*{x}}_0,\vb{P}_0$};
    \node (stop) [startstop, below of=prop] {Stop};
    
    \draw [arrow] (start) -- (input);
    \draw [arrow] (input) -- (iod);
    \draw [arrow] (iod) -- (proj);
    \draw [arrow] (proj) -- (prune);
    \draw [arrow] (prune) -- (outlier);
    \draw [arrow] (outlier) -- node[anchor=east] {no} (merge);
    \draw [arrow] (outlier) -| node[anchor=south] {yes} (restore);
    \draw[arrow] (merge) -- (loop);
    \draw [arrow] (restore) |- (loop);
    \draw [arrow] (loop) -|  node[anchor=north] {no} ++ (2.5cm,0) |- (repl);
    \draw [arrow] (outlier) -- node[anchor=east] {no} (merge);
    \draw[arrow] (loop) -- node[anchor=east] {yes} (prop);
    \draw [arrow] (prop) -| (count);
    \draw[arrow] (count) |- (proj);
    \draw [arrow] (repl) -- (extract);
    \draw [arrow] (extract) -- (od);
    \draw [arrow] (od) -- (out);
    \draw [arrow] (out) |- ++ (0,-2cm) -| (stop);
    
    \end{tikzpicture}
    \caption{Workflow of the proposed robust \glsentryshort{od} solution}
    \label{fig:flowchart}
\end{figure}
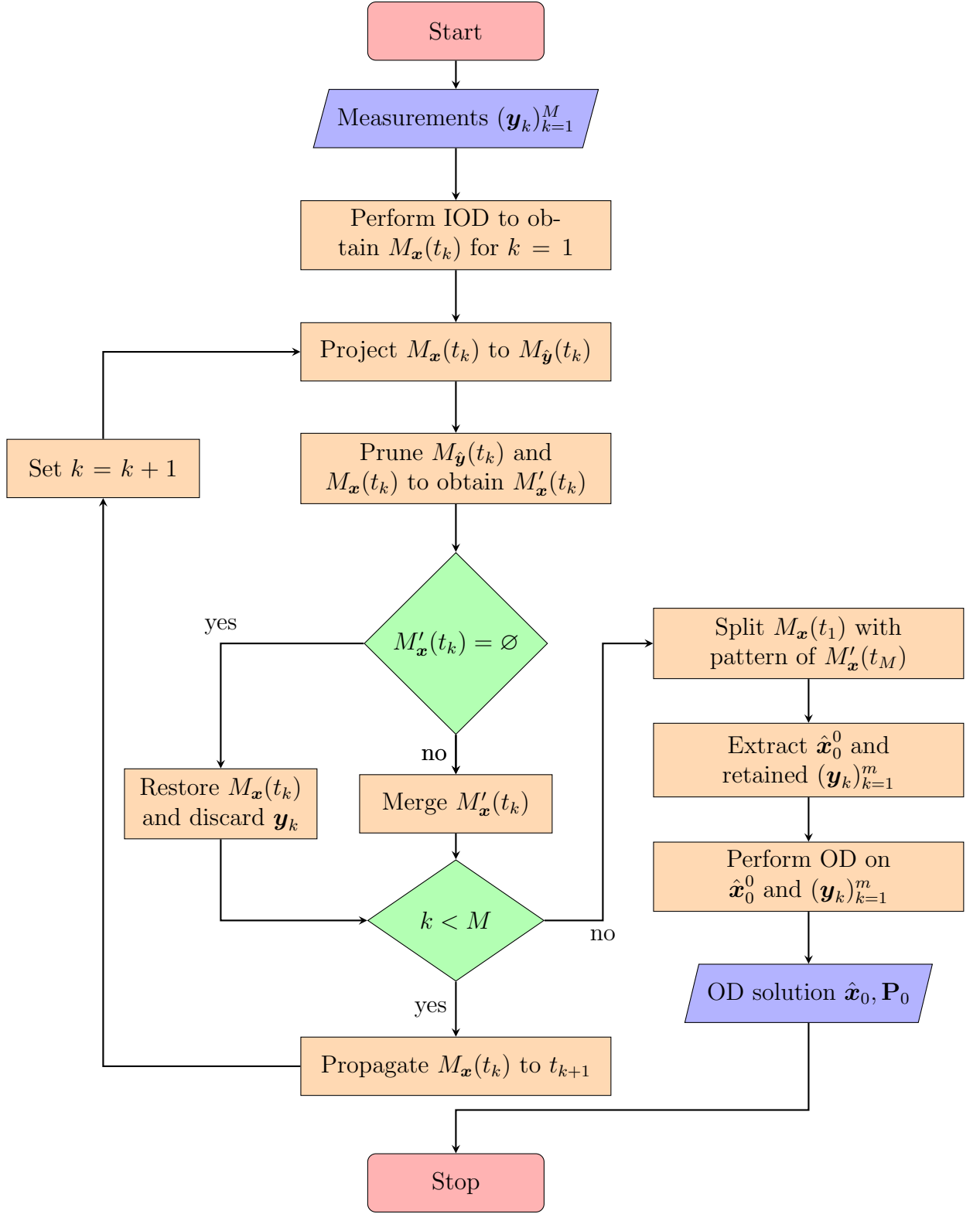

At convergence, the solution to the \glsentryshort{od} problem consists in an estimated state $\hat{\vb*{x}}_0$ that minimizes the cost function \cref{eq:cost_fun_ls} or \cref{eq:cost_fun_lsar} and an associated covariance matrix $\vb{P}_0$. If a new set of observations become available, it is then conceivable to repeat the process and compute an updated estimate $\hat{\vb*{x}}'_0$ and covariance $\vb*{P}'_0$ by taking into account the newly available information. Rather than repeating each step in \cref{fig:flowchart} by considering the two batches of measurements at once, it is more efficient to initialize $M_{\vb*{x}}(t_0)$ from $\hat{\vb*{x}}_0,\vb*{P}_0$ and run the pruning and \glsentryshort{od} algorithms only on the second batch. Assuming a Gaussian distribution for $\hat{\vb*{x}}_0$, the last is initialized as the \glsentryshort{da} vector
\begin{equation}
\begin{aligned}
    [\vb{x}_0] &= \hat{\vb*{x}}_0 + \vb*{V}_0\cdot\left(c\sqrt{\vb*{\lambda}_0}\delta\vb*{x}_0 \right)\\
    &= \mathcal{T}_{\vb{x}_0}(\delta\vb*{x}_0)
\end{aligned}
\label{eq:init_poly_second_batch}
\end{equation}
with $c$ the $z$-score, $\vb*{P}_0=\vb*{V}_0\vb*{\Lambda}_0\vb*{V}_0^T$ the covariance matrix eigendecomposition and $\vb*{\lambda}_0$ the main diagonal of $\vb*{\Lambda}_0$. In this case $M_{\vb*{x}}(t_0)$ will contain a single element, i.e. $[\vb{x}_0]$, which is a multivariate polynomial in the six independent \glsentryshort{da} variables $\delta\vb*{x}_0$. It is thus different from \cref{sec:seq_algo} where the initial manifold at $t_0$ is function of the angular variables $(\delta\vb*{\alpha},\delta\vb*{\delta})$ introduced in \cref{sec:iod} and might contain more than one element if splits occurred during \glsentryshort{iod}. Moreover, representing the state estimate with a mean vector $\hat{\vb*{x}}_0$ and a covariance matrix $\vb*{P}_0$ is a limitation of the current algorithm, since all second order information gained with the \glsentryshort{da}-based \glsentryshort{iod} and pruning routines is lost when the batch \glsentryshort{od} algorithm is started. Future developments will thus focus on replacing the current \glsentryshort{od} routines with high-order \glsentryshort{da}-based sequential filters~\cite{Servadio2020}. The pruning and \glsentryshort{od} steps could then be combined in a single sequential algorithm in which the state update is also performed in the algebra of Taylor polynomials. The final output would then be a second-order manifold $M_{\vb*{x}}(t_M)$, and newly available observations for $t>t_M$ would be seamlessly processed without the needs for a re-initialization of the \glsentryshort{da} variables as in \cref{eq:init_poly_second_batch}.

\section{Numerical Application}\label{sec:application}

This section presents an application of the developed algorithm in real operational scenarios. The selected object is the rocket body of an Ariane 44L launched in 2002 and catalogued with \glsentryshort{norad} ID 27402. This object is on a \glsentryfull{gto} with semi-major axis of about \SI{23000}{\km} and inclination of \SI{3.9}{\degree}, and is regularly tracked by the optical telescopes of the \glsentryshort{tarot} network~\cite{Boer2017}. This scenario was selected since \glsentryfullpl{heo} are notoriously challenging for both \glsentryshort{iod} and \glsentryshort{up} tasks. Specifically for the object of interest, observations are available only in the proximity of the apogee where the slow motion of the \glsentryshort{so}, coupled with a large distance from the sensors, leads to small angular separations between measurements and thus large uncertainties in the computed \glsentryshort{iod} solution. \Glsentryshortpl{gto} are then characterized by very fast perigee passages, during which the objects is subject to strong perturbations from the Earth's atmosphere, and slow apogee passages, during which the effects of atmospheric drag are negligible, while \glsentrylong{srp} and lunisolar perturbations prevail. This continuous alternation of fast and slow dynamics, coupled with a large eccentricity that magnifies its nonlinearities, make more difficult to maintain an accurate representation of the state over time. A successful application of the algorithm to this scenario will thus give confidence on its robustness and applicability to a wide spectrum of orbital regimes.

Firstly, the performance of \glsentryshort{da}-based \glsentryshort{up} methods is assessed in \cref{sec:mf_validation} by comparing their accuracy and computational efficiency against a \glsentrylong{hf} \glsentryshort{mc} simulation. Two analyses are then conducted in \cref{sec:synthetic_measurements,sec:tarot_measurements} to highlight different aspects of the proposed tool. In both cases, two objects are considered: the one of interest (the target), and a third-party one (the outlier). The outlier is simulated, and noised measurements are produced out of it to simulate observations that should not be associated to the target object. Instead, the observations that should be associated to the \glsentryshort{so} of interest are either generated from a synthetic trajectory in \cref{sec:synthetic_measurements} or drawn from real observations in \cref{sec:tarot_measurements}. Both analyses consider a time span of about five days that includes five passages of the \glsentryshortpl{so} over the two telescopes in La Reunion and Calern (France). Passages epochs, duration and number of measurements are summarized in \cref{tab:passages}.

\begin{table}[H]
    \centering
    \caption{Epoch, duration and number of measurements for the considered passages}
    \begin{tabular}{ccccc}
    \toprule
    Sensor & Epoch (UTC) & Duration & \#measurements\\ \midrule
    La Reunion & 2019-02-25T18:49:01.148 & 05:13:10.189 & 8\\
    La Reunion & 2019-02-27T23:27:20.511 & 00:00:48.972 & 3\\
    La Reunion & 2019-03-01T01:35:55.561 & 00:00:48.019 & 3\\
    La Reunion & 2019-03-01T23:10:59.277 & 00:00:24.011 & 2\\
    Calern & 2019-03-01T23:58:27.006 & 00:00:24.926 & 2\\ \bottomrule
    \end{tabular}
    \label{tab:passages}
\end{table}

In the followings the analytical \glsentryshort{sgp} propagator~\cite{Vallado2006a} is used as \glsentrylong{lf} dynamical model. Instead, a \glsentrylong{hf} numerical propagator is set up with the force models listed below

\begin{itemize}
    \setlength{\itemsep}{0pt plus 1pt}
    \item Earth's gravity potential modeled with spherical harmonics truncated at order/degree 16~\cite{Holmes2002}
    \item Sun and Moon third-body attraction with bodies' position from \glsentryshort{jpl} \glsentryshort{de}440 ephemerides~\cite{Park2021}
    \item Atmospheric drag with Harris-Priester atmosphere model~\cite{Harris1962}
    \item \glsentrylong{srp} with cylindrical Earth's shadow model~\cite{Hubaux2012}
\end{itemize}

Taylor polynomials in \cref{eq:iod_alpha_delta_da,eq:pn_cov_da,eq:real_meas_da} are initialized with a $z$-score equal to $c=3$ such that the $99.73\%$ of the probability mass is accurately mapped through each nonlinear transformation using the \glsentryshort{loads} scheme. The algorithms were coded in Java and make use of the following external libraries
\begin{itemize}
    \item \glsentryshort{cnes}'s \glsentryfullinv{pace} library as \glsentryshort{da} backend to perform operations on Taylor polynomials
    \item Orekit\footnote{\href{https://www.orekit.org/}{https://www.orekit.org/}} for standard flight dynamics routines (\glsentryshort{sgp} propagator, force models, time scales and reference frames conversions)
    \item Hipparchus\footnote{\href{https://hipparchus.org/}{https://hipparchus.org/}} for the 8(5,3) Dormand-Prince numerical integrator and the \glsentrylong{lm} optimizer
    \item MOSEK\footnote{\href{https://www.mosek.com/}{https://www.mosek.com/}} interfaced through its Fusion API for Java\footnote{\href{https://docs.mosek.com/latest/javafusion/index.html}{https://docs.mosek.com/latest/javafusion/index.html}} as \glsentryshort{lpp} solver
\end{itemize}
\noindent
All simulations were run on a cluster equipped with Skylake Intel\textsuperscript{\textregistered} Xeon\textsuperscript{\textregistered} Gold 6126 @ \SI{2.6}{\giga\hertz} and \SI{96}{\giga\byte} of RAM running Red Hat\textsuperscript{\textregistered} Enterprise Linux\textsuperscript{\textregistered} 7.8.

\subsection{Validation of multifidelity method}\label{sec:mf_validation}

The performance of the \glsentrylong{mf} \glsentryshort{up} method presented in \cref{sec:mf_up} is discussed in this section to justify its adoption within the developed \glsentryshort{od} framework. Four different \glsentryshort{up} schemes are considered, namely
\begin{enumerate}
    \item \Glsentryfull{lf} \glsentryshort{da}-based propagation (first step of \cref{sec:mf_up})
    \item \Glsentryfull{mf} \glsentryshort{da}-based propagation (as used within the \glsentryshort{od} framework)
    \item \Glsentryfull{hf} \glsentryshort{da}-based propagation (\glsentryshort{da}-aware numerical propagator directly embedded in the \glsentryshort{loads} framework)
    \item \glsentryshort{mc} simulation with $N_s=10^4$ samples using the \glsentrylong{hf} dynamics
\end{enumerate}
A comparison is then performed in terms of both accuracy and computational efficiency, with the \glsentryshort{mc} simulation taken as ground truth. A manifold of \glsentryshortpl{ic} $M(t_0)=\{[\vb*{x}^{(i)}(t_0)]\}$ is firstly obtained by solving the \glsentryshort{iod} problem as described in \cref{sec:iod} using the real measurements from the first passage in \cref{tab:passages}. These \glsentryshortpl{ic} are then propagated over a five days time span with the four methods listed above. First, random samples $\{\vb*{x}^{(l)}(t_0)\}$ for $l=1,\ldots,N_s$ are drawn from $M(t_0)$ and propagated \glsentrylong{pw} to obtain the reference empirical distribution $\{\vb*{x}^{(l)}(t_f)\}$ at final time. Then, $M(t_0)$ is propagated in the \glsentryshort{da}/\glsentryshort{loads} framework to obtain the three manifolds $M_{LF}(t_f),M_{MF}(t_f)$ and $M_{HF}(t_f)$ corresponding to the three \glsentryshort{da}-based \glsentryshort{up} schemes. These manifolds are then evaluated on $\{\vb*{x}^{(l)}(t_0)\}$ to obtain $\{\vb*{x}_{LF}^{(l)}(t_f)\},\{\vb*{x}_{MF}^{(l)}(t_f)\}$ and $\{\vb*{x}_{HF}^{(l)}(t_f)\}$ without the needs for further propagation. The \glsentryfull{rmse} is finally used to assess the accuracy of the polynomial-based \glsentryshort{up} methods. The last is defined as
\begin{equation}
    \vb*{e}_{RMSE}=\sqrt{\dfrac{1}{N_s}\sum_{l=1}^{N_s}\left(\vb*{x}^{(l)}-\hat{\vb*{x}}^{(l)}\right)^2}
    \label{eqn:rmse_def}
\end{equation}
with $\{\vb*{x}^{(l)}\}$ the expected samples, i.e. the output of the \glsentryshort{mc} simulation, and $\{\hat{\vb*{x}}^{(l)}\}$ the actual samples, i.e. the result of the polynomial evaluation of either of the three manifolds at final time.

\begin{table}[!ht]
    \centering
    \captionsetup{justification=centering}
    \caption{\Glsentrylong{rmse} for a five days propagation}
    \begin{adjustwidth}{-2.5cm}{-2.5cm}
    \centering
    \small
    \sisetup{round-mode=places,round-precision=3}
    \setlength\tabcolsep{2pt}
    \begin{tabular}{r *6{S[scientific-notation=true,table-format=1.3e-1]}}
    \toprule
    & {$x$, \si{\km}} & {$y$, \si{\km}} & {$z$, \si{\km}} & {$v_x$, \si{\km\per\s}} & {$v_y$, \si{\km\per\s}} & {$v_z$, \si{\km\per\s}}\\ \midrule
    \glsentryshort{lf} & 181.501633403219 & 121.951277883348 & 6.50414537496815 & 0.0319013130988913 & 0.0105220416669801 & 0.000184099443676873\\
    \glsentryshort{mf} & 1.81079564950228 & 1.32140575553321 & 0.0855189236968127 & 0.000301435053817822 & 0.000101002771932697 & 2.81198772476723e-06\\
    \glsentryshort{hf} & 0.140434694853268 & 0.101660903665224 & 0.00679731568611594 & 2.2908733910041e-05 & 7.56303651096141e-06 & 4.01360588736488e-07\\
    \bottomrule
    \end{tabular}
    \end{adjustwidth}
    \label{tab:rmse_prop}
\end{table}
The \glsentryshortpl{rmse} for the case under consideration are reported in \cref{tab:rmse_prop}. As seen from this table, the \glsentrylong{lf} dynamics is unable to accurately capture the time evolution of the state distribution, resulting in average errors that exceed \SI{218}{\km} in position and \SI{3.36e-2}{\km\per\s} in velocity. The results are greatly improved with the introduction of the \glsentrylong{mf} correction, with \glsentryshortpl{rmse} of about \SI{2.24}{\km} and \SI{3.18e-4}{\km\per\s} respectively. The last are further reduced with the \glsentrylong{hf} \glsentryshort{da}-based propagation whose accuracy is in the order of \SI{1.74e-1}{\km} in position and \SI{2.41e-5}{\km\per\s} in velocity, an order of magnitude better than the \glsentrylong{mf} solution.
\begin{table}[!ht]
    \centering
    \captionsetup{justification=centering}
    \caption{Computational time for a five days propagation}
    \begin{adjustwidth}{-2.5cm}{-2.5cm}
    \centering
    \small
    \sisetup{round-mode=places,round-precision=3}
    \setlength\tabcolsep{2pt}
    \begin{tabular}{r *5{S[scientific-notation=false,table-format=8.3]}}
    \toprule
    & {\Glsentrylong{lf}} & {\Glsentrylong{mf}} & {\Glsentrylong{hf}} & {Single sample} & {\Glsentrylong{mc}}\\ \midrule
    $t$, \si{\s} & 0.237 & 3.392 & 144.934 & 3.0409875 & 30409.875\\
    $t/t_{MF}$ & 0.0698702830188679 & 1 & 42.7281839622642 & 0.896517541273585 & 8965.17541273585\\
    \bottomrule
    \end{tabular}
    \end{adjustwidth}
    \label{tab:runtime_prop}
\end{table}
However, as shown in \cref{tab:runtime_prop}, additional accuracy comes at a high computational cost. The \glsentrylong{hf} propagator is in fact 43 times slower than its \glsentrylong{mf} counterpart, which is in turn 14 times slower than the \glsentrylong{lf} one. Recalling that the last coincides with the first step of the \glsentrylong{mf} scheme, it is possible to conclude that the computational cost of the hybrid method is driven by the correction step, i.e. the \glsentrylong{pw} propagation of the polynomials' constant parts in \glsentrylong{hf} dynamics. The overall cost of the \glsentrylong{mf} scheme is thus comparable to that of a \glsentrylong{hf} \glsentryshort{mc} simulation with a number of samples equal to the number of domains generated by the \glsentryshort{loads} algorithm. In this case no splits are triggered during the propagation in alternate equinoctial elements, and a single split is required to map the uncertainty onto Cartesian space. The runtime of the hybrid scheme is thus similar to that of the \glsentryshort{lf} simulation plus the time required to propagate a single \glsentryshort{mc} sample in \glsentrylong{hf}. This fact can be easily verified by comparing the first, second and fourth columns of \cref{tab:runtime_prop}.

\begin{figure}[!ht]
    \centering
    \includegraphics[width=0.6\textwidth]{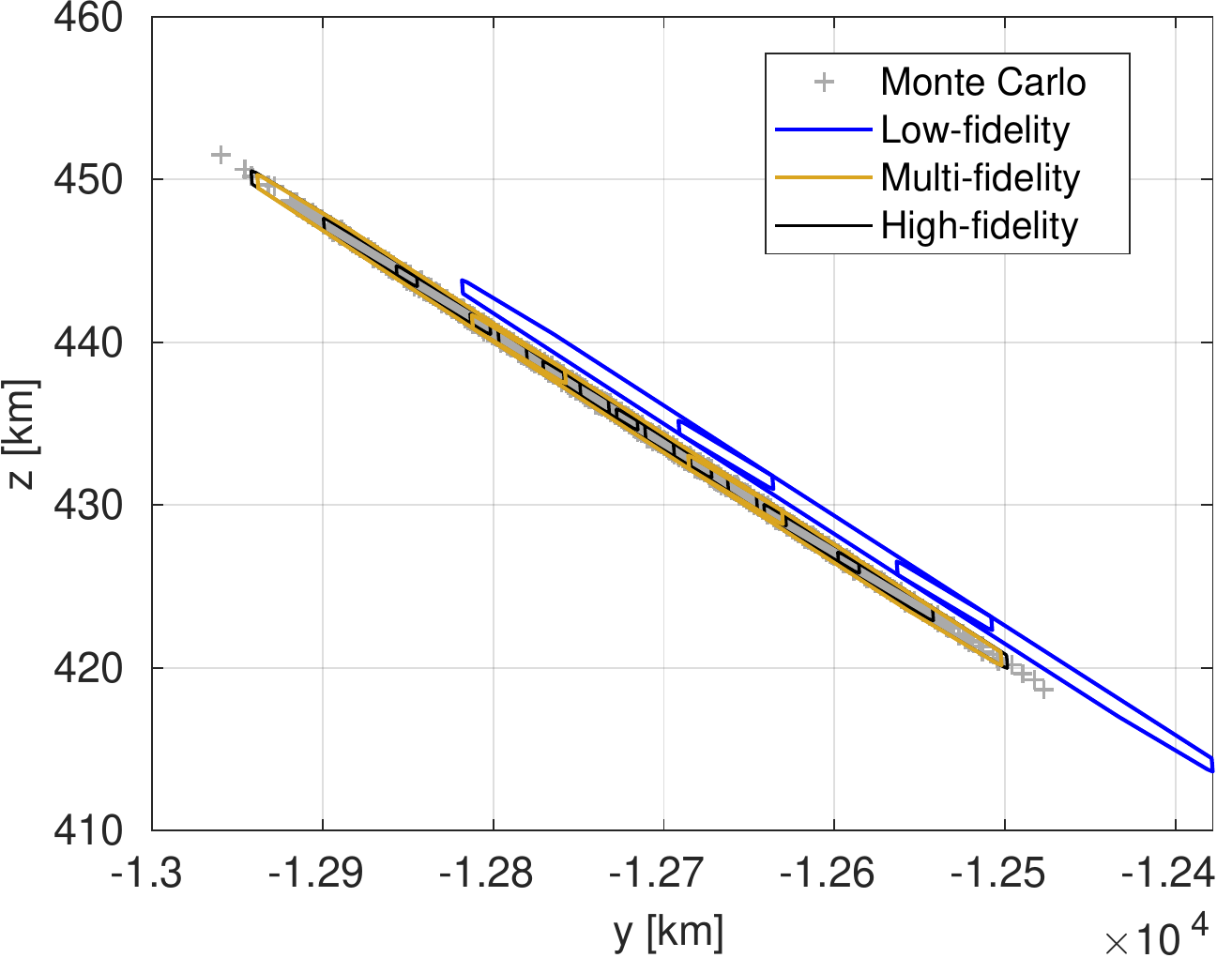}
    \captionsetup{justification=centering}
    \caption{Projection onto the $y-z$ plane of the \glsentrylong{mc} samples and the polynomial domains edges after a five days propagation}
    \label{fig:samp_edge_prop_y_z}
\end{figure}
\begin{figure}[!ht]
    \centering
    \includegraphics[width=0.6\textwidth]{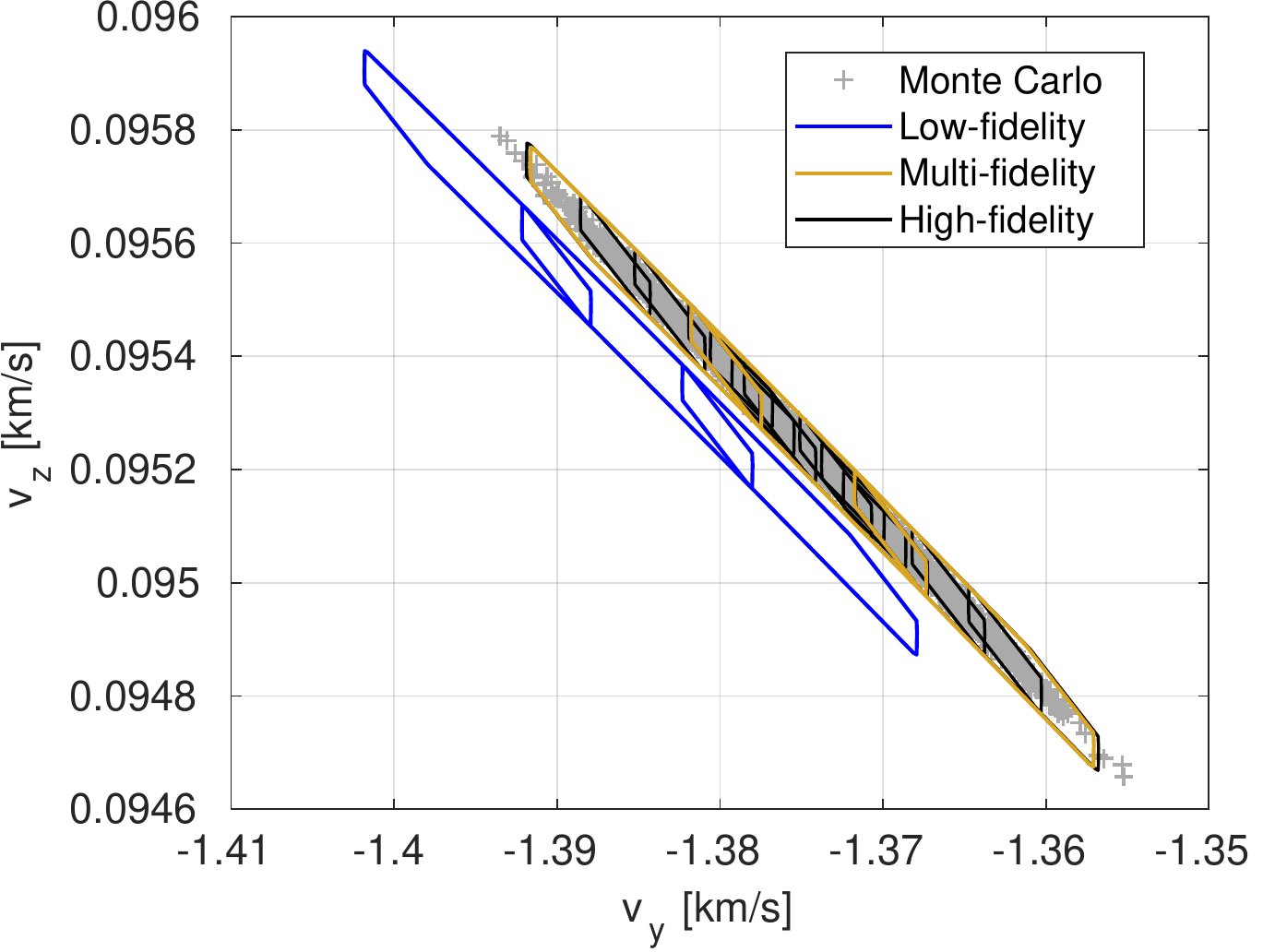}
    \captionsetup{justification=centering}
    \caption{Projection onto the $v_y-v_z$ plane of the \glsentrylong{mc} samples and the polynomial domains edges after a five days propagation}
    \label{fig:samp_edge_prop_vy_vz}
\end{figure}

\Cref{fig:samp_edge_prop_y_z,fig:samp_edge_prop_vy_vz} depicts the \glsentryshort{mc} samples (gray plus signs) and the edges of the polynomial domains for the four cases enumerated before. Domains in blue correspond to the \glsentrylong{lf} solution, which is clearly shifted with respect to the ground truth in both $y-z$ and $v_y-v_z$ projections. Instead, domains in gold (\glsentrylong{mf} solution) and in black (\glsentrylong{hf} solution) overlap almost perfectly and provide a visual demonstration of the accuracy of the proposed \glsentrylong{mf} propagation algorithm. As expected, most of the \glsentryshort{mc} samples fall inside these edges thus demonstrating the validity of the \glsentryshort{da}-based \glsentryshort{up} methods. Moreover, the good accuracy of the \glsentrylong{mf} method can be explained by examining the different steps involved to obtain the final solution. The first one consists in propagating the initial uncertainty in \glsentrylong{lf}, and the output is equivalent to the \glsentrylong{lf} solution depicted in blue. As seen in both pictures, the union of the regions mapped by these polynomials has the same shape of that mapped by the \glsentrylong{hf} ones (delimited by black edges), and the only noticeable difference is a shift between the two. It is thus possible to conclude that \glsentryshort{sgp}, even though it fails in capturing the evolution of the absolute dynamics, is capable of accurately capture the relative motion in the neighborhoods of the polynomials' constant parts delimited by the depicted edges. The second step is instead a propagation of the polynomials' constant part in \glsentrylong{hf} using the same dynamical model as the \glsentryshort{mc} simulation. Once these two solutions are available, the \glsentrylong{hf} trajectories are taken as reference with respect to which the uncertainty is expressed, and a simple substitution of the polynomials' centers allows to obtain a solution that matches the accuracy of the \glsentrylong{hf} one. To conclude, the \glsentrylong{mf} solution provides the best trade-off for the use case under investigation. The small loss in accuracy observed in \cref{tab:rmse_prop} is in fact compensated by the large computational gain shown in \cref{tab:runtime_prop}.

\subsection{Synthetic measurements}\label{sec:synthetic_measurements}

The first analysis of the complete \glsentryshort{od} chain uses synthetic measurements for both objects, so as a true reference trajectory is available to assess the algorithm's performance. For the target \glsentryshort{so} an assumed true state is retrieved from \glsentryshort{tle} data. The outlier orbit is then obtained by perturbing the target's eccentricity by $\Delta e=-0.02$. The Keplerian parameters at epoch are given in \cref{tab:simu_ics} for both \glsentryshortpl{so}. The trajectories of the two objects are simulated using the \glsentrylong{hf} propagator and synthetic measurements are drawn from them. These are then corrupted with a zero-mean Gaussian white noise with standard deviations $\sigma_{\alpha}=\SI{1.285}{\arcsecond}$ and $\sigma_{\delta}=\SI{1.280}{\arcsecond}$ to simulate sensor noise. Four different scenarios are then analyzed as listed in \cref{tab:sim_scenarios}. In scenario A and B only correlated measurements are processed. These are meant to validate the tools before further complexity is added to the problem. Observations associated to the outlier are then added in scenarios C and D.

\begin{table}[!ht]
    \centering
    \caption{\Glsentrylongpl{ic} for target and outlier orbits}
    \sisetup{round-mode=places,round-precision=3}
    \begin{tabular}{c S[scientific-notation=false,table-format=5.3] *2{S[scientific-notation=false,table-format=1.3]} S[scientific-notation=false,table-format=3.3] S[scientific-notation=false,table-format=-3.3] S[scientific-notation=false,table-format=2.3]}
    \toprule
    & {$a$, \si{\km}} & {$e$, -} & {$i$, \si{\degree}} & {$\omega$, \si{\degree}} & {$\Omega$, \si{\degree}} & {$M_0$, \si{\degree}}\\ \midrule
    Target & 22953.852669768778 & 0.707854612716 & 3.387521317683 & 172.980213527756 & -168.891499315499 & 60.742995057860\\
    Outlier & 22953.852669768778 & 0.687854612716 & 3.387521317683 & 172.980213527756 & -168.891499315499 & 60.742995057860\\
    \bottomrule
    \end{tabular}
    \label{tab:simu_ics}
\end{table}

\begin{table}[!ht]
    \centering
    \caption{Test case scenarios with simulated measurements}
    \begin{tabular}{ccc}
    \toprule
    Scenario \# & Pruning & Outliers\\ \midrule
    A & NO & NO\\
    B & YES & NO\\
    C & NO & 3\textsuperscript{rd} passage\\
    D & YES & 3\textsuperscript{rd} passage\\ \bottomrule
    \end{tabular}
    \label{tab:sim_scenarios}
\end{table}

\begin{figure}[!ht]
    \centering
    \includegraphics[width=0.6\textwidth]{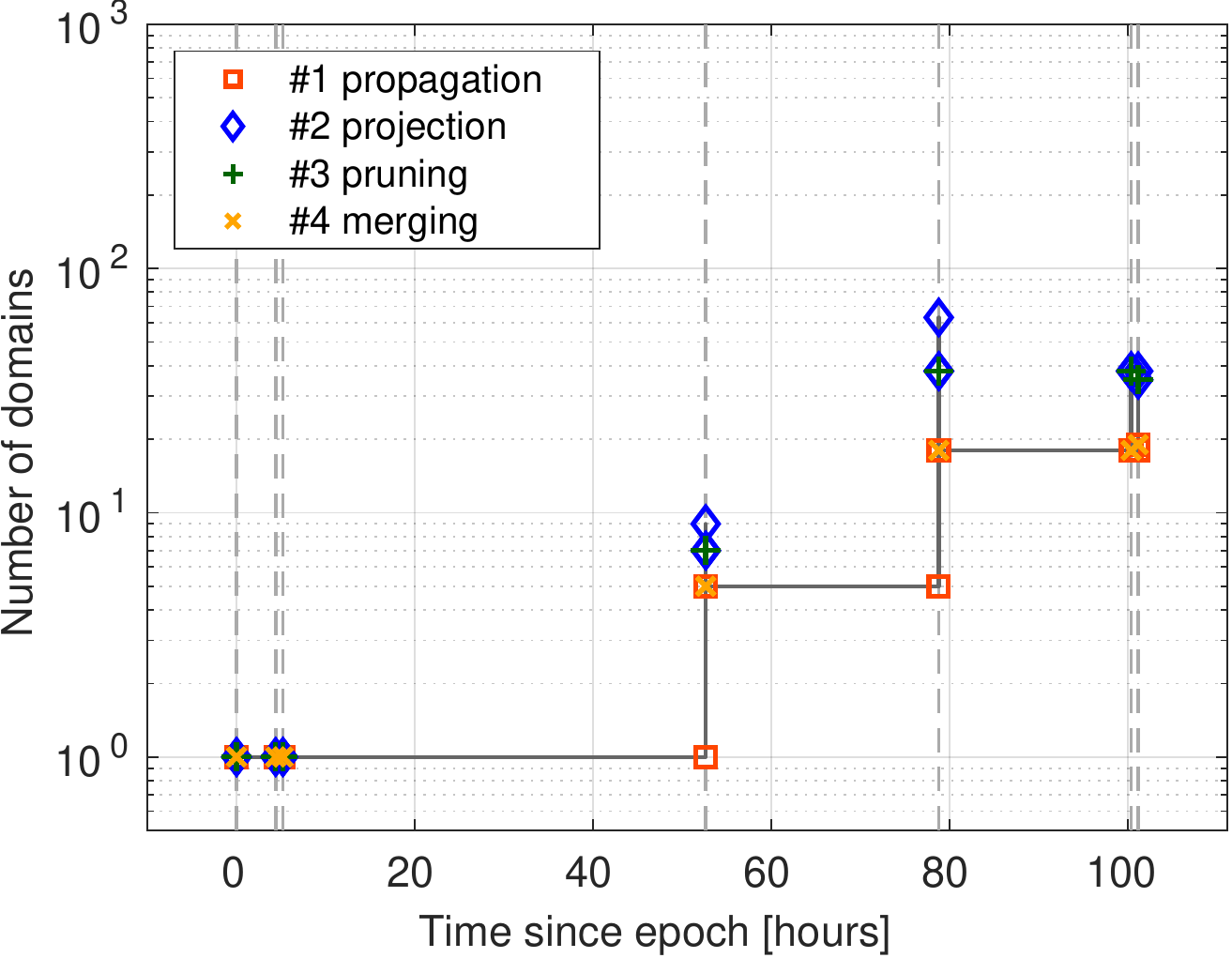}
    \caption{Number of domains as function of propagation time for scenario B}
    \label{fig:hist_simu_no_out}
\end{figure}

\begin{figure}[!ht]
    \centering
    \includegraphics[width=0.6\textwidth]{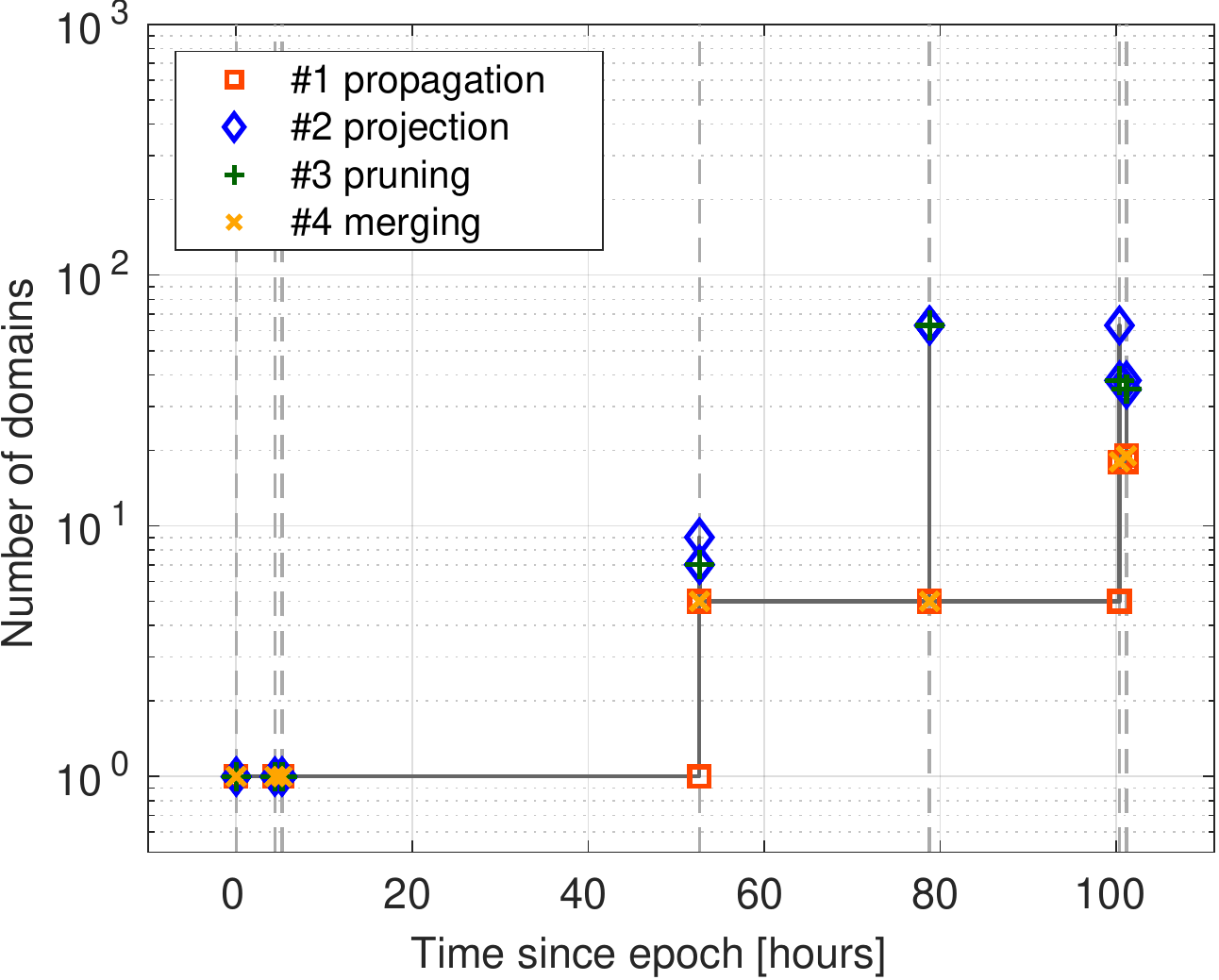}
    \caption{Number of domains as function of propagation time for scenario D}
    \label{fig:hist_simu_out}
\end{figure}

The \glsentryshort{iod} algorithm described in \cref{sec:iod} is run on the first passage and the estimated state at $t_0$ used to initialize the pruning scheme. In this case, the \glsentryshort{loads} algorithm generates a single domain and thus $M_{\vb*{x}}=\{[\vb*{x}_0]\}$. In scenarios B and D this manifold is then propagated to the latest epoch using the algorithm in \cref{sec:seq_algo}. The evolution over time of the number of domains is summarized in \cref{tab:hist_simu_no_out,tab:hist_simu_out} respectively. A visual representation is also given in \cref{fig:hist_simu_no_out,fig:hist_simu_out}. In this pictures, orange squares and yellow crosses correspond to the size of the osculating elements manifolds after propagation (step \#1) and merging (step \#4). Blue diamonds and green plus signs correspond to the size of the observables manifolds after projection (step \#2) and pruning (step \#3). Vertical dashed lines correspond to the observations epochs.

\begin{table}[!ht]
    \centering
    \caption{Number of domains as function of propagation time for scenario B}
    \small
    \sisetup{round-mode=places,round-precision=3}
    \begin{tabular}{S[scientific-notation=false,table-format=3.3] cccc}
    \toprule
    & \multicolumn{4}{c}{Number of domains}\\ \midrule
    {Time since $t_0$, \si{\hour}} & \#1 Propagation & \#2 Projection & \#3 Pruning & \#4 Merging\\ \midrule
    0 & 1 & 1 & 1 & 1\\
    0.00667388888888889 & 1 & 1 & 1 & 1\\
    0.0133258333333333 & 1 & 1 & 1 & 1\\
    4.39337888888889 & 1 & 1 & 1 & 1\\
    4.40004833333333 & 1 & 1 & 1 & 1\\
    5.20617527777778 & 1 & 1 & 1 & 1\\
    5.21283222222222 & 1 & 1 & 1 & 1\\
    5.21949694444444 & 1 & 1 & 1 & 1\\ \midrule
    52.6387119444444 & 1 & 9 & 7 & 5\\
    52.6453727777778 & 5 & 7 & 7 & 5\\
    52.6523152777778 & 5 & 7 & 7 & 5\\ \midrule
    78.7817813888889 & 5 & 63 & 38 & 18\\
    78.7884552777778 & 18 & 38 & 38 & 18\\
    78.79512 & 18 & 38 & 38 & 18\\ \midrule
    100.366146944444 & 18 & 38 & 38 & 18\\
    100.372816666667 & 18 & 38 & 38 & 18\\ \midrule
    101.157182777778 & 18 & 38 & 35 & 19\\
    101.164106666667 & 19 & 35 & 35 & 19\\ \bottomrule
    \end{tabular}
    \label{tab:hist_simu_no_out}
\end{table}

\begin{table}[!ht]
    \centering
    \caption{Number of domains as function of propagation time for scenario D}
    \small
    \sisetup{round-mode=places,round-precision=3}
    \begin{tabular}{S[scientific-notation=false,table-format=3.3] cccc}
    \toprule
    & \multicolumn{4}{c}{Number of domains}\\ \midrule
    {Time since $t_0$, \si{\hour}} & \#1 Propagation & \#2 Projection & \#3 Pruning & \#4 Merging\\ \midrule
    0 & 1 & 1 & 1 & 1\\
    0.00667388888888889 & 1 & 1 & 1 & 1\\
    0.0133258333333333 & 1 & 1 & 1 & 1\\
    4.39337888888889 & 1 & 1 & 1 & 1\\
    4.40004833333333 & 1 & 1 & 1 & 1\\
    5.20617527777778 & 1 & 1 & 1 & 1\\
    5.21283222222222 & 1 & 1 & 1 & 1\\
    5.21949694444444 & 1 & 1 & 1 & 1\\ \midrule
    52.6387119444444 & 1 & 9 & 7 & 5\\
    52.6453727777778 & 5 & 7 & 7 & 5\\
    52.6523152777778 & 5 & 7 & 7 & 5\\ \midrule
    78.7817813888889 & 5 & 63 & 63 & 5\\
    78.7884552777778 & 5 & 63 & 63 & 5\\
    78.79512 & 5 & 63 & 63 & 5\\ \midrule
    100.366146944444 & 5 & 63 & 38 & 18\\
    100.372816666667 & 18 & 38 & 38 & 18\\ \midrule
    101.157182777778 & 18 & 38 & 35 & 19\\
    101.164106666667 & 19 & 35 & 35 & 19\\ \bottomrule
    \end{tabular}
    \label{tab:hist_simu_out}
\end{table}

The effectiveness of the pruning algorithm is clearly seen by comparing, for the same epoch, the number of domains after steps \#2 and \#3. It is also noted that the merging step is key in limiting the number of domains processed in step \#1. If not performed, the propagation algorithm will enter the next iteration with a number of domains corresponding to the third column (or green plus signs) rather than the fourth one (or yellow crosses). The effectiveness of the selected state representation in limiting the number of domains generated during propagation is also demonstrated. All splits are in fact triggered by the projection onto the observables space and the increased number of domains to be propagated in subsequent steps is the result of an imperfect merging during step \#4. Since domains are discarded in step \#3, several incomplete triplets are found by the merging scheme which thus fails to recombine the polynomials as they were before step \#2. The correct detection of outliers is then demonstrated in \cref{tab:hist_simu_out}, where no modification of the domains occur in the third measurement block which corresponds to uncorrelated measurements. In this case, the algorithm correctly detects the outliers and does not prune the projected manifold.

The success rate of the pruning algorithm in correctly identifying outlier measurements is reported in \cref{tab:pruning_succes_simu} for scenarios without outliers (cases A and B) and with outliers (cases C and D) respectively. This table shows the number of correctly associated (true positives) and discarded (true negative) measurements as well as the number of those that should have been marked as correlated (false negatives) or outliers (false positives) but they have been not. Results obtained with both \glsentrylong{lf} and \glsentrylong{mf} \glsentryshort{up} methods are shown. The \glsentrylong{hf} solution is instead omitted since a perfect match is already obtained with the \glsentrylong{mf} one. The same is however not true for the \glsentrylong{lf} method for which five false negatives are observed in both scenarios. The importance of the \glsentrylong{mf} correction is thus once again remarked by these data.
\begin{table}[!ht]
    \centering
    \caption{Success rate of outlier detection with simulated measurements}
    \begin{adjustwidth}{-2.5cm}{-2.5cm}
    \small
    \centering
    \begin{tabular}{cccccc}
    \toprule
    Case \# & Dynamics & True positives & False negatives & False positives & True negatives\\ \midrule
    \multirow{2}{*}{A,B} & \glsentryshort{lf} & 13 & 5 & 0 & 0\\
    & \glsentryshort{mf} & 18 & 0 & 0 & 0\\
    \midrule
    \multirow{2}{*}{C,D} & \glsentryshort{lf} & 10 & 5 & 0 & 3\\
    & \glsentryshort{mf} & 15 & 0 & 0 & 3\\
    \bottomrule
    \end{tabular}
    \end{adjustwidth}
    \label{tab:pruning_succes_simu}
\end{table}

\begin{figure}[!ht]
    \centering
    \includegraphics[width=0.6\textwidth]{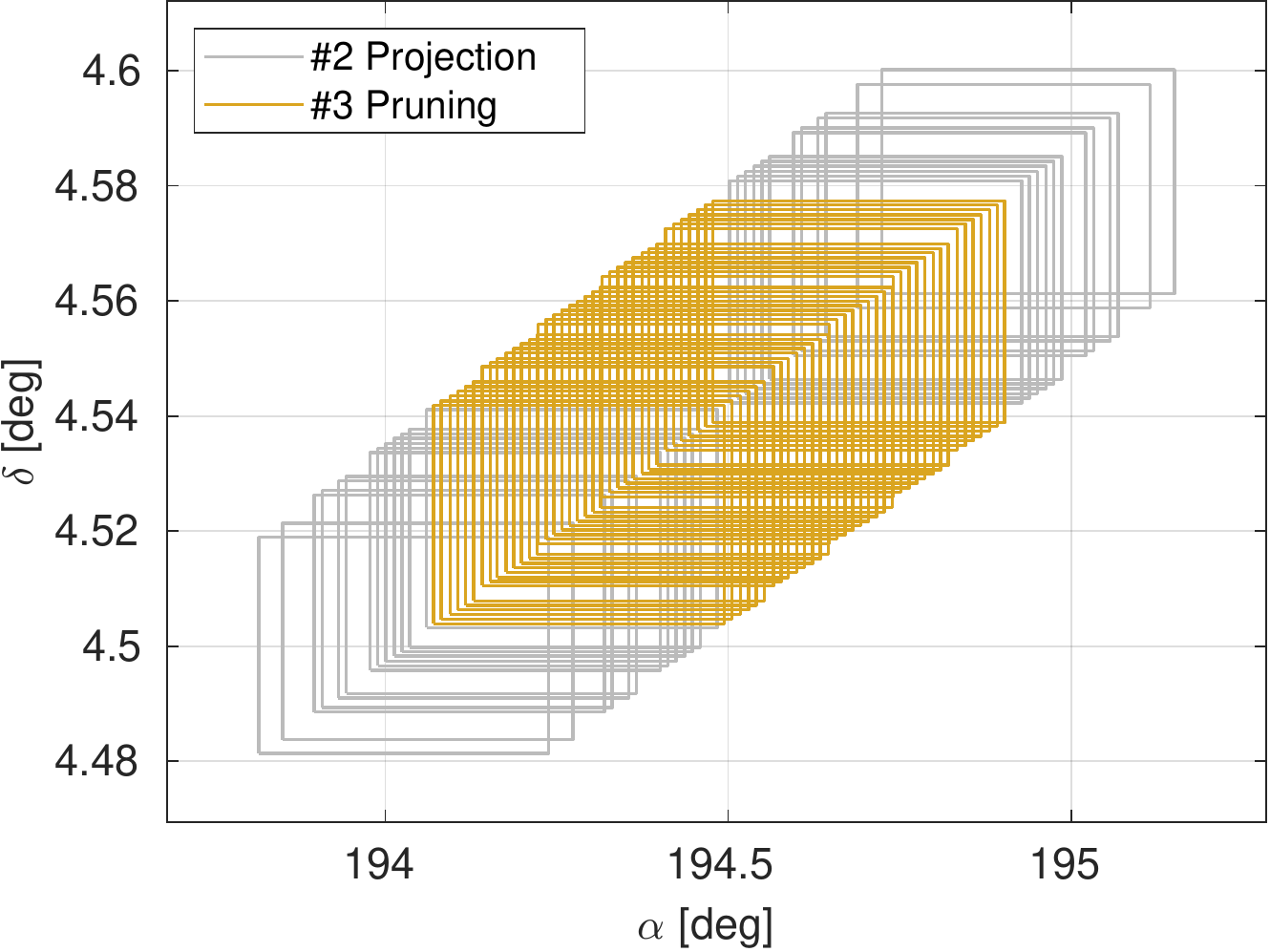}
    \captionsetup{justification=centering}
    \caption{Estimated bounds of projected and retained domains at \SI{100.366}{\hour} after $t_0$ for scenario D}
    \label{fig:bounds_simu_out}
\end{figure}

\begin{figure}[!ht]
    \centering
    \includegraphics[width=0.6\textwidth]{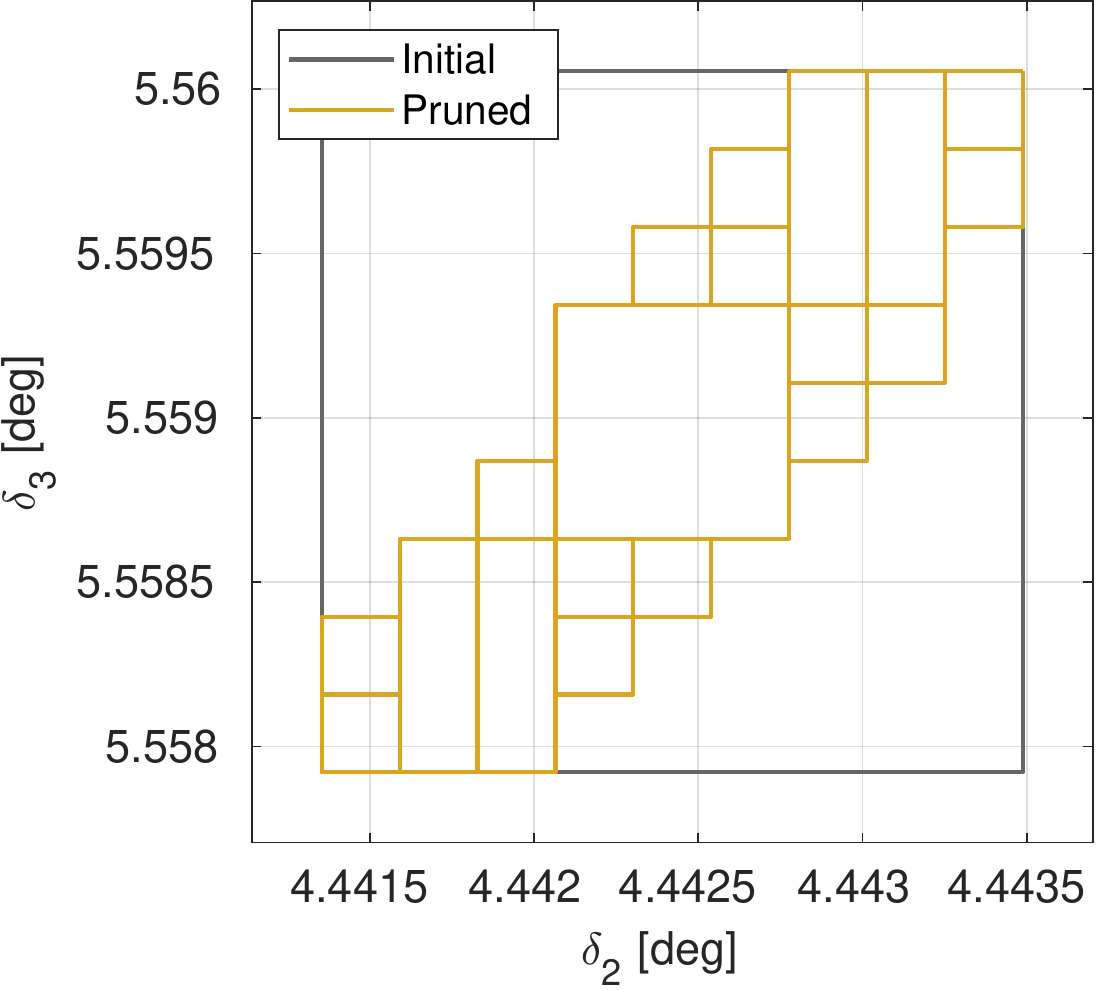}
    \captionsetup{justification=centering}
    \caption{Edges of the base \glsentryshort{iod} solution and retained domains in $(\delta_2,\delta_3)$ space for scenario D}
    \label{fig:iod_simu_out}
\end{figure}

An example of pruning is presented in \cref{fig:bounds_simu_out} where the estimated bounds for both $M_{\hat{\vb*{y}}}(t_k)$ (gray boxes) and $M_{\hat{\vb*{y}}}'(t_k)$ (gold boxes) are plotted in $(\alpha,\delta)$ space. The large overlap between gold boxes is an indicator that most of the uncertainty is on the \glsentryshort{los} distance which cannot be observed with optical telescopes. Still, the pruning algorithm is highly effective in reducing the uncertainty on the observed quantities. This is demonstrated also in \cref{fig:iod_simu_out} which represents the uncertainty on the \glsentryshort{iod} solution in $(\delta_2,\delta_3)$ space. The gray box identifies the edges of the initial domain which corresponds to the uncertainty on the measurements. Instead, the golden boxes are the edges of the domains that maps into the solution manifold at the last epoch. It is thus seen the effectiveness of the algorithm in tightening the bounds on the initial estimate as the propagation steps forward and more information is exploited to refine the \glsentryshort{iod} solution.

\Glsentryshort{od} is then performed for all scenarios comparing the performance of the \glsentryshort{ls} and \glsentryshort{lsar} algorithms introduced in \cref{sec:batch_od}. The true Cartesian state of the target object and the \glsentryshort{iod} solution used to initialize the \glsentryshort{od} routines are given in \cref{tab:od_true_simu}. The resulting estimation errors and $\pm 3\sigma$ bounds are summarized in \cref{tab:od_sol_simu_no_out,tab:od_sol_simu_out_no_prune,tab:od_sol_simu_out_prune} for scenarios A and B, C and D respectively.

\begin{table}[!ht]
    \centering
    \captionsetup{justification=centering}
    \caption{True state and initial guess for scenarios A to D}
    \begin{adjustwidth}{-2.5cm}{-2.5cm}
    \centering
    \small
    \sisetup{round-mode=places,round-precision=5}
    \setlength\tabcolsep{2pt}
    \begin{tabular}{r *6{S[scientific-notation=true,table-format=-1.5e-1]}}
    \toprule
    & {$x$, \si{\km}} & {$y$, \si{\km}} & {$z$, \si{\km}} & {$v_x$, \si{\km\per\s}} & {$v_y$, \si{\km\per\s}} & {$v_z$, \si{\km\per\s}}\\ \midrule
    True & -21551.184664630193 & 14404.866452074804 & -1082.462558770526 & -3.580403901491 & -0.736464589895 & 0.001943794765\\
    Guess & -21553.629554466817 & 14406.082502818883 & -1082.362628418208 & -3.579945419480 & -0.736567402660 & 0.001921790734\\ \bottomrule
    \end{tabular}
    \end{adjustwidth}
    \label{tab:od_true_simu}
\end{table}

\Cref{tab:od_sol_simu_no_out} corresponds to the two scenarios with no outliers in which the pruning algorithm has no effect on the final solution. Both algorithms provide a consistent estimate of the \glsentryshort{so} state converging in three and two iterations for the \glsentryshort{ls} and \glsentryshort{lsar} methods respectively. As expected, the \glsentryshort{ls} solution is characterized by tighter bounds on the states since the algorithm is the minimum variance estimator for a model with zero-mean Gaussian noise distribution~\cite{Prabhu2021}.

\begin{table}[!ht]
    \centering
    \captionsetup{justification=centering}
    \caption{\glsentryshort{ls} and \glsentryshort{lsar} estimation errors and $\pm 3\sigma$ bounds for scenario A and B}
    \begin{adjustwidth}{-2.5cm}{-2.5cm}
    \centering
    \small
    \sisetup{round-mode=places,round-precision=3}
    \setlength\tabcolsep{2pt}
    \begin{tabular}{r *6{S[scientific-notation=true,table-format=1.3e-1]}}
    \toprule
    & {$x$, \si{\km}} & {$y$, \si{\km}} & {$z$, \si{\km}} & {$v_x$, \si{\km\per\s}} & {$v_y$, \si{\km\per\s}} & {$v_z$, \si{\km\per\s}}\\ \midrule
    \glsentryshort{ls} Error & 0.00143288266633762 & 0.166254060672003 & 0.0640439894502724 & 1.43662669525573e-05 & 6.95245707222816e-06 & 8.22763419768502e-07\\
    $\pm 3\sigma$ bounds & 0.911943829770987 & 0.456451793548531 & 0.205293174909563 & 0.000171046554685692 & 4.02257709877014e-05 & 1.53497227064365e-05\\ \midrule
    \glsentryshort{lsar} Error & 0.180380981372442 & 0.082141679550848 & 0.045135145490439 & 1.90461080920885e-05 & 1.05962850008136e-05 & 2.08297951491271e-07\\
    $\pm 3\sigma$ bounds & 1.13040992399803 & 0.565739521827316 & 0.254563548634605 & 0.000212009774878624 & 4.98719155696088e-05 & 1.90331893420941e-05\\ \bottomrule
    \end{tabular}
    \end{adjustwidth}
    \label{tab:od_sol_simu_no_out}
\end{table}

The effects of measurement outliers on the estimated state is then shown in \cref{tab:od_sol_simu_out_no_prune}. With no measurement preprocessing, the \glsentryshort{ls} estimator fails to converge to the true solution and returns an estimate characterized by large errors in both position ($\approx\SI{700}{\km}$) and velocity ($\approx\SI{0.12}{\km\per\s}$) converging in six iterations. At the same time, the variance on the estimated state is close to that of scenarios A and B, thus placing the true state outside the $\pm 3\sigma$ bounds. On the contrary, the \glsentryshort{lsar} algorithm is capable of rejecting the outliers and still converging to the true solution at the expenses of a longer runtime (five iterations vs two in the previous cases). The estimated bounds are also consistent to those obtained in scenarios A and B and they thus include the true state reported in \cref{tab:od_true_simu}. 

\begin{table}[!ht]
    \centering
    \captionsetup{justification=centering}
    \caption{\glsentryshort{ls} and \glsentryshort{lsar} estimation errors and $\pm 3\sigma$ bounds for scenario C}
    \begin{adjustwidth}{-2.5cm}{-2.5cm}
    \centering
    \small
    \sisetup{round-mode=places,round-precision=3}
    \setlength\tabcolsep{2pt}
    \begin{tabular}{r *6{S[scientific-notation=true,table-format=1.3e-1]}}
    \toprule
    & {$x$, \si{\km}} & {$y$, \si{\km}} & {$z$, \si{\km}} & {$v_x$, \si{\km\per\s}} & {$v_y$, \si{\km\per\s}} & {$v_z$, \si{\km\per\s}}\\ \midrule
    \glsentryshort{ls} Error & 579.080210298213 & 377.552242082303 & 34.4826671195127 & 0.118092221182537 & 0.0246243820305996 & 0.00256862784270661\\
    $\pm 3\sigma$ bounds & 0.861941148124937 & 0.422408271524913 & 0.210029303800098 & 0.000157272736644034 & 3.68025114179516e-05 & 1.47256660988291e-05\\ \midrule
    \glsentryshort{lsar} Error & 0.772982128250064 & 0.257594372595664 & 0.0148850772888534 & 0.000138442689863824 & 4.15625713590147e-05 & 3.4100442515069e-07\\
    $\pm 3\sigma$ bounds & 1.13015101083691 & 0.565569011378571 & 0.254564917720344 & 0.000211949556622479 & 4.98638470585985e-05 & 1.903236267704e-05\\ \bottomrule
    \end{tabular}
    \end{adjustwidth}
    \label{tab:od_sol_simu_out_no_prune}
\end{table}

\begin{table}[!ht]
    \centering
    \captionsetup{justification=centering}
    \caption{\glsentryshort{ls} and \glsentryshort{lsar} estimation errors and $\pm 3\sigma$ bounds for scenario D}
    \begin{adjustwidth}{-2.5cm}{-2.5cm}
    \centering
    \small
    \sisetup{round-mode=places,round-precision=3}
    \setlength\tabcolsep{2pt}
    \begin{tabular}{r *6{S[scientific-notation=true,table-format=-1.3e-1]}}
    \toprule
    & {$x$, \si{\km}} & {$y$, \si{\km}} & {$z$, \si{\km}} & {$v_x$, \si{\km\per\s}} & {$v_y$, \si{\km\per\s}} & {$v_z$, \si{\km\per\s}}\\ \midrule
    \glsentryshort{ls} Error & 0.188526166259416 & 0.043656216609869 & 0.0440300238728799 & 2.44510472610899e-05 & 1.51214609080539e-05 & 2.04168313632785e-06\\
    $\pm 3\sigma$ bounds & 1.09632762203042 & 0.598195189390871 & 0.226457257104369 & 0.000211436830082596 & 4.79484577763055e-05 & 1.67898624395595e-05\\ \midrule
    \glsentryshort{lsar} Error & 0.333753328664473 & 0.00578430230075149 & 0.0555757121099707 & 4.97227651961683e-05 & 1.85939889899553e-05 & 4.3010131765688e-07\\
    $\pm 3\sigma$ bounds & 1.35934111088477 & 0.741705821914885 & 0.280806887552135 & 0.000262159299521155 & 5.94544810458999e-05 & 2.08191449226171e-05\\ \bottomrule
    \end{tabular}
    \end{adjustwidth}
    \label{tab:od_sol_simu_out_prune}
\end{table}

\begin{table}[!ht]
    \centering
    \captionsetup{justification=centering}
    \caption{Runtime and number of iterations with simulated measurements}
    \begin{adjustwidth}{-2.5cm}{-2.5cm}
    \centering
    \small
    \sisetup{round-mode=places,round-precision=3}
    \setlength\tabcolsep{2pt}
    \begin{tabular}{c *3{S[scientific-notation=false,table-format=3.3]} cc}
    \toprule
    & \multicolumn{3}{c}{$t$, \si{\s}} & \multicolumn{2}{c}{\# iterations}\\ \midrule
    Case \# & Pruning & LS & LSAR & LS & LSAR\\ \midrule
    A & {-} & 331.642 & 125.777 & 3 & 2\\
    B & 70.297 & 331.257 & 125.944 & 3 & 2\\
    C & {-} & 642.552 & 324.101 & 6 & 5\\
    D & 40.718 & 133.349 & 87.541 & 2 & 2\\ \bottomrule
    \end{tabular}
    \end{adjustwidth}
    \label{tab:od_runtime_simu}
\end{table}

The accuracy of the \glsentryshort{ls} estimate is then restored in scenario D in which outliers have been detected and discarded by the pruning scheme. Results in \cref{tab:od_sol_simu_out_prune} are in fact similar to those in \cref{tab:od_sol_simu_no_out} and the \glsentryshort{ls} estimation is again characterized by tighter bounds on the states. The \glsentryshort{lsar} algorithm also benefits from measurements preprocessing since excluding the uncorrelated observations from the processed ones reduces the number of iterations from five to just two. This results in a total runtime of about 40\% of that of scenario C as can be verified by inspecting in \cref{tab:od_runtime_simu}.

\subsection{\texorpdfstring{\glsentryshort{tarot}}{TAROT} measurements}\label{sec:tarot_measurements}

The second analysis uses real observation data from the \glsentryshort{tarot} network to estimate the target \glsentryshort{so} state. Synthetic measurements for the outlier orbit generated in \cref{sec:synthetic_measurements} are still used to corrupt the available data and assess the robustness of the proposed tool. The four scenarios presented in \cref{sec:synthetic_measurements} are reproduced here as summarized in \cref{tab:true_scenarios}. To account for the unknown forces acting on the \glsentryshortpl{so},  stochastic accelerations modeled as additive white Gaussian noise are included in the \glsentrylong{hf} dynamics. Their effect is taken into account by the pruning scheme using the \glsentryshort{snc} algorithm described in \cref{sec:snc}. A deterministic dynamics is instead used for the subsequent batch \glsentryshort{od}.

\begin{table}[!ht]
    \centering
    \caption{Test case scenarios with real measurements from the \glsentryshort{tarot} network}
    \begin{tabular}{ccc}
    \toprule
    Case \# & Pruning & Outliers\\ \midrule
    E & NO & NO\\
    F & YES & NO\\
    G & NO & 3\textsuperscript{rd} passage\\
    H & YES & 3\textsuperscript{rd} passage\\ \bottomrule
    \end{tabular}
    \label{tab:true_scenarios}
\end{table}

\begin{figure}[!ht]
    \centering
    \includegraphics[width=0.6\textwidth]{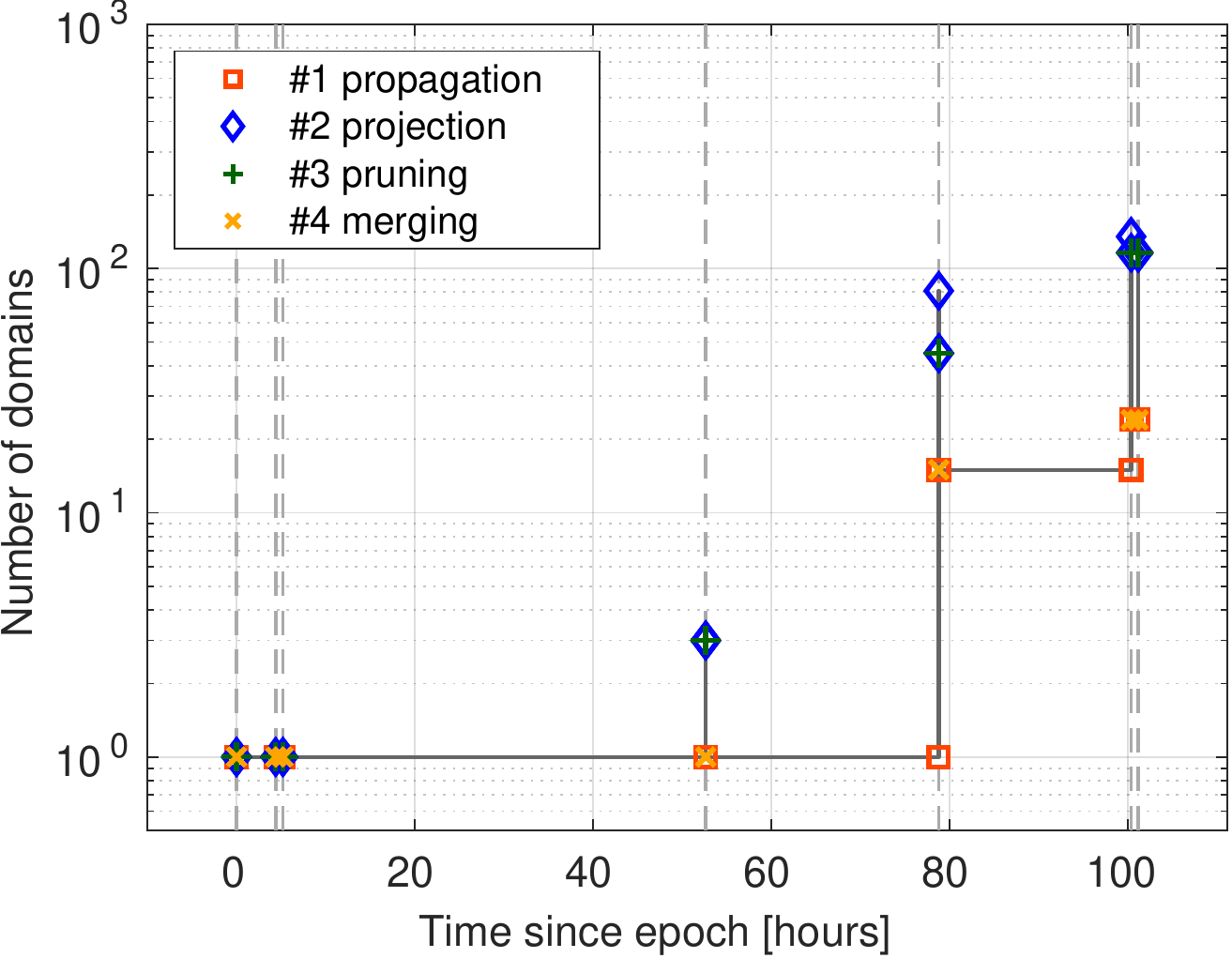}
    \caption{Number of domains as function of propagation time for scenario F}
    \label{fig:hist_true_no_out}
\end{figure}

\begin{figure}[!ht]
    \centering
    \includegraphics[width=0.6\textwidth]{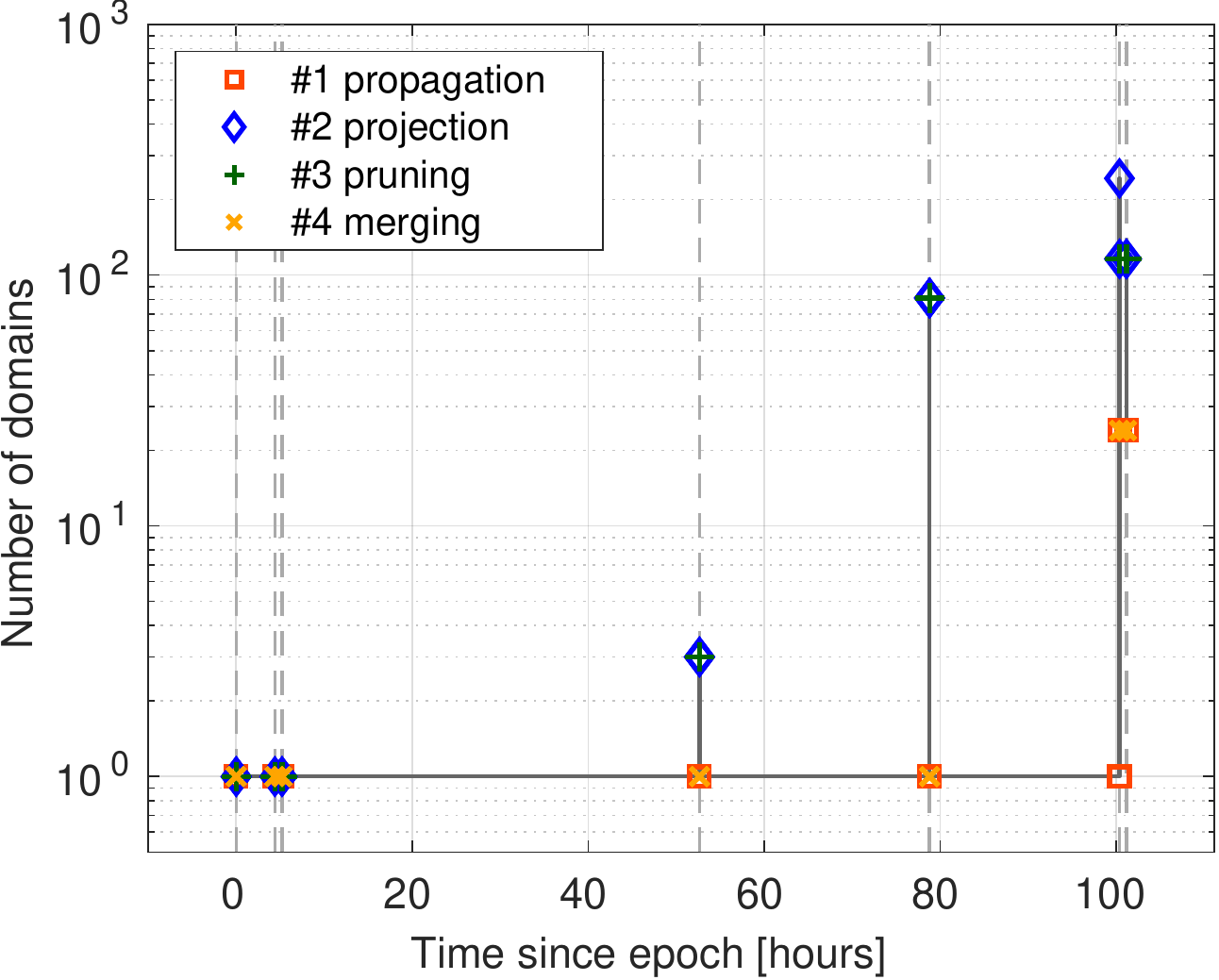}
    \caption{Number of domains as function of propagation time for scenario H}
    \label{fig:hist_true_out}
\end{figure}

The \glsentryshort{iod} algorithm run on the first passage still return a single domain which is then used to initialize the pruning scheme in scenarios F and H. The time evolution of the number of domains is summarized in \cref{tab:hist_true_no_out,tab:hist_true_out} and depicted in \cref{fig:hist_true_no_out,fig:hist_true_out} respectively for the two scenarios. The immediate effect of stochastic accelerations can be identified in the larger number of domains generated during the projection (\#2) step. Indeed, integrating the uncertainty due to these perturbations results in an inflated covariance for the prior states which, in turns, triggers more splits when projected onto the observables space. At the same time, the pruning scheme is even more effective and the ratio between projected and retained domains is higher in \cref{tab:hist_true_no_out} than in \cref{tab:hist_simu_no_out}. The success rate of the pruning algorithm in correctly identify outlier measurements is summarized in \cref{tab:pruning_succes_tarot} and is identical to that of \cref{tab:pruning_succes_simu}. The effectiveness of the \glsentrylong{mf} \glsentryshort{up} method in providing an accurate estimate of the prior at a fraction of the computational cost of a \glsentrylong{hf} propagation is thus demonstrated also in this case.

\begin{table}[!ht]
    \centering
    \caption{Number of domains as function of propagation time for scenario F}
    \small
    \sisetup{round-mode=places,round-precision=3}
    \begin{tabular}{S[scientific-notation=false,table-format=3.3] cccc}
    \toprule
    & \multicolumn{4}{c}{Number of domains}\\ \midrule
    {Time since $t_0$, \si{\hour}} & \#1 Propagation & \#2 Projection & \#3 Pruning & \#4 Merging\\ \midrule
    0 & 1 & 1 & 1 & 1\\
    0.00667388888888889 & 1 & 1 & 1 & 1\\
    0.0133258333333333 & 1 & 1 & 1 & 1\\
    4.39337888888889 & 1 & 1 & 1 & 1\\
    4.40004833333333 & 1 & 1 & 1 & 1\\
    5.20617527777778 & 1 & 1 & 1 & 1\\
    5.21283222222222 & 1 & 1 & 1 & 1\\
    5.21949694444444 & 1 & 1 & 1 & 1\\ \midrule
    52.6387119444444 & 1 & 3 & 3 & 1\\
    52.6453727777778 & 1 & 3 & 3 & 1\\
    52.6523152777778 & 1 & 3 & 3 & 1\\ \midrule
    78.7817813888889 & 1 & 81 & 45 & 15\\
    78.7884552777778 & 15 & 45 & 45 & 15\\
    78.79512 & 15 & 45 & 45 & 15\\ \midrule
    100.366146944444 & 15 & 135 & 116 & 24\\
    100.372816666667 & 24 & 116 & 116 & 24\\ \midrule
    101.157182777778 & 24 & 116 & 116 & 24\\
    101.164106666667 & 24 & 116 & 116 & 24\\ \bottomrule
    \end{tabular}
    \label{tab:hist_true_no_out}
\end{table}

\begin{table}[!ht]
    \centering
    \caption{Number of domains as function of propagation time for scenario H}
    \small
    \sisetup{round-mode=places,round-precision=3}
    \begin{tabular}{S[scientific-notation=false,table-format=3.3] cccc}
    \toprule
    & \multicolumn{4}{c}{Number of domains}\\ \midrule
    {Time since $t_0$, \si{\hour}} & \#1 Propagation & \#2 Projection & \#3 Pruning & \#4 Merging\\ \midrule
    0 & 1 & 1 & 1 & 1\\
    0.00667388888888889 & 1 & 1 & 1 & 1\\
    0.0133258333333333 & 1 & 1 & 1 & 1\\
    4.39337888888889 & 1 & 1 & 1 & 1\\
    4.40004833333333 & 1 & 1 & 1 & 1\\
    5.20617527777778 & 1 & 1 & 1 & 1\\
    5.21283222222222 & 1 & 1 & 1 & 1\\
    5.21949694444444 & 1 & 1 & 1 & 1\\ \midrule
    52.6387119444444 & 1 & 3 & 3 & 1\\
    52.6453727777778 & 1 & 3 & 3 & 1\\
    52.6523152777778 & 1 & 3 & 3 & 1\\ \midrule
    78.7817813888889 & 1 & 81 & 81 & 1\\
    78.7884552777778 & 1 & 81 & 81 & 1\\
    78.79512 & 1 & 81 & 81 & 1\\ \midrule
    100.366146944444 & 1 & 243 & 116 & 24\\
    100.372816666667 & 24 & 116 & 116 & 24\\ \midrule
    101.157182777778 & 24 & 116 & 116 & 24\\
    101.164106666667 & 24 & 116 & 116 & 24\\ \bottomrule

    \end{tabular}
    \label{tab:hist_true_out}
\end{table}

\begin{table}[!ht]
    \centering
    \caption{Success rate of outlier detection with real measurements}
    \begin{adjustwidth}{-2.5cm}{-2.5cm}
    \small
    \centering
    \begin{tabular}{cccccc}
    \toprule
    Case \# & Dynamics & True positives & False negatives & False positives & True negatives\\ \midrule
    \multirow{2}{*}{E,F} & \glsentryshort{lf} & 13 & 5 & 0 & 0\\
    & \glsentryshort{mf} & 18 & 0 & 0 & 0\\
    \midrule
    \multirow{2}{*}{G,H} & \glsentryshort{lf} & 10 & 5 & 0 & 3\\
    & \glsentryshort{mf} & 15 & 0 & 0 & 3\\
    \bottomrule
    \end{tabular}
    \end{adjustwidth}
    \label{tab:pruning_succes_tarot}
\end{table}

\begin{figure}[!ht]
    \centering
    \includegraphics[width=0.6\textwidth]{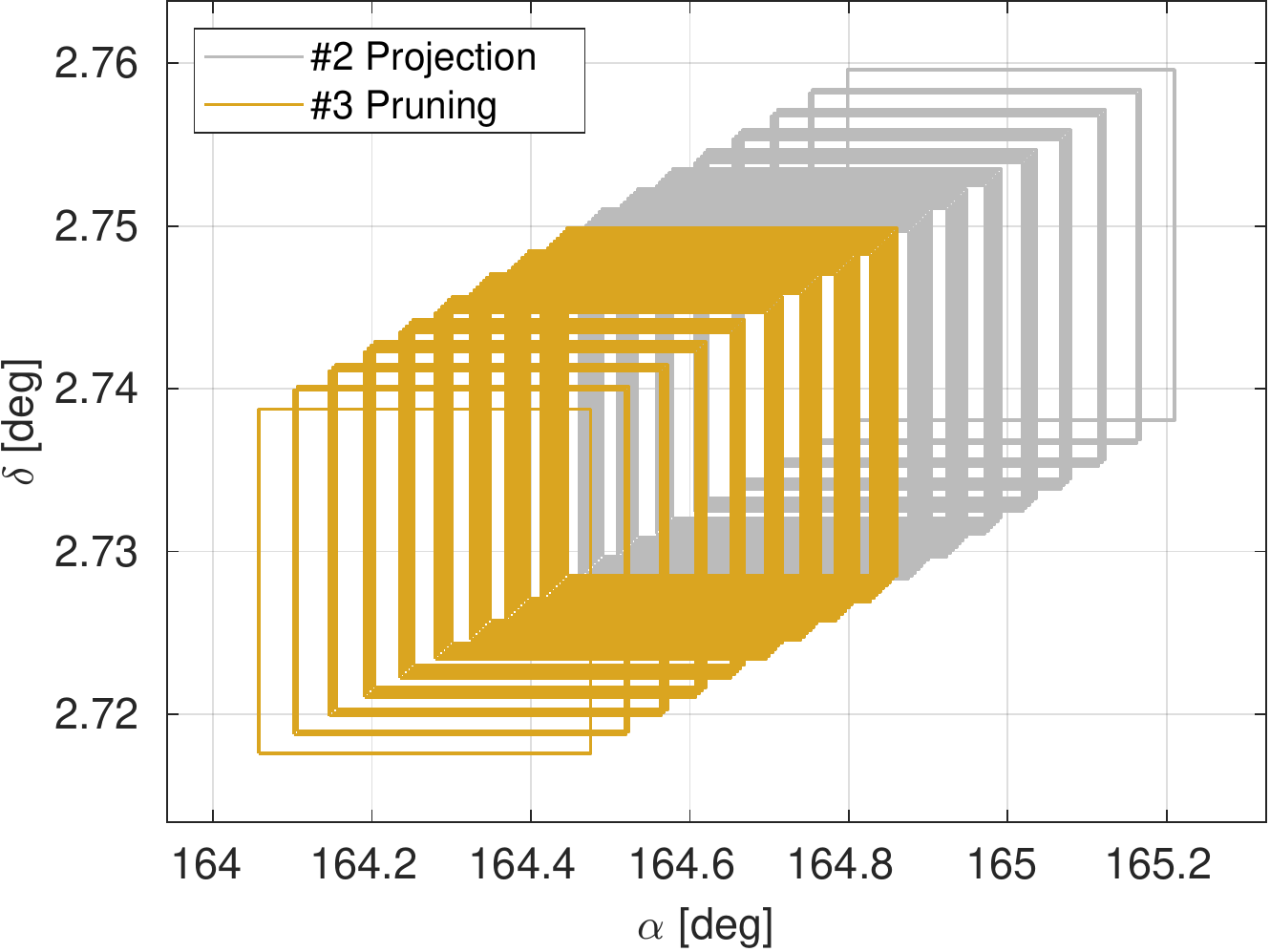}
    \captionsetup{justification=centering}
    \caption{Estimated bounds of projected and retained domains at \SI{78.782}{\hour} after $t_0$ for scenario F}
    \label{fig:bounds_true_no_out}
\end{figure}

\begin{figure}[!ht]
    \centering
    \includegraphics[width=0.6\textwidth]{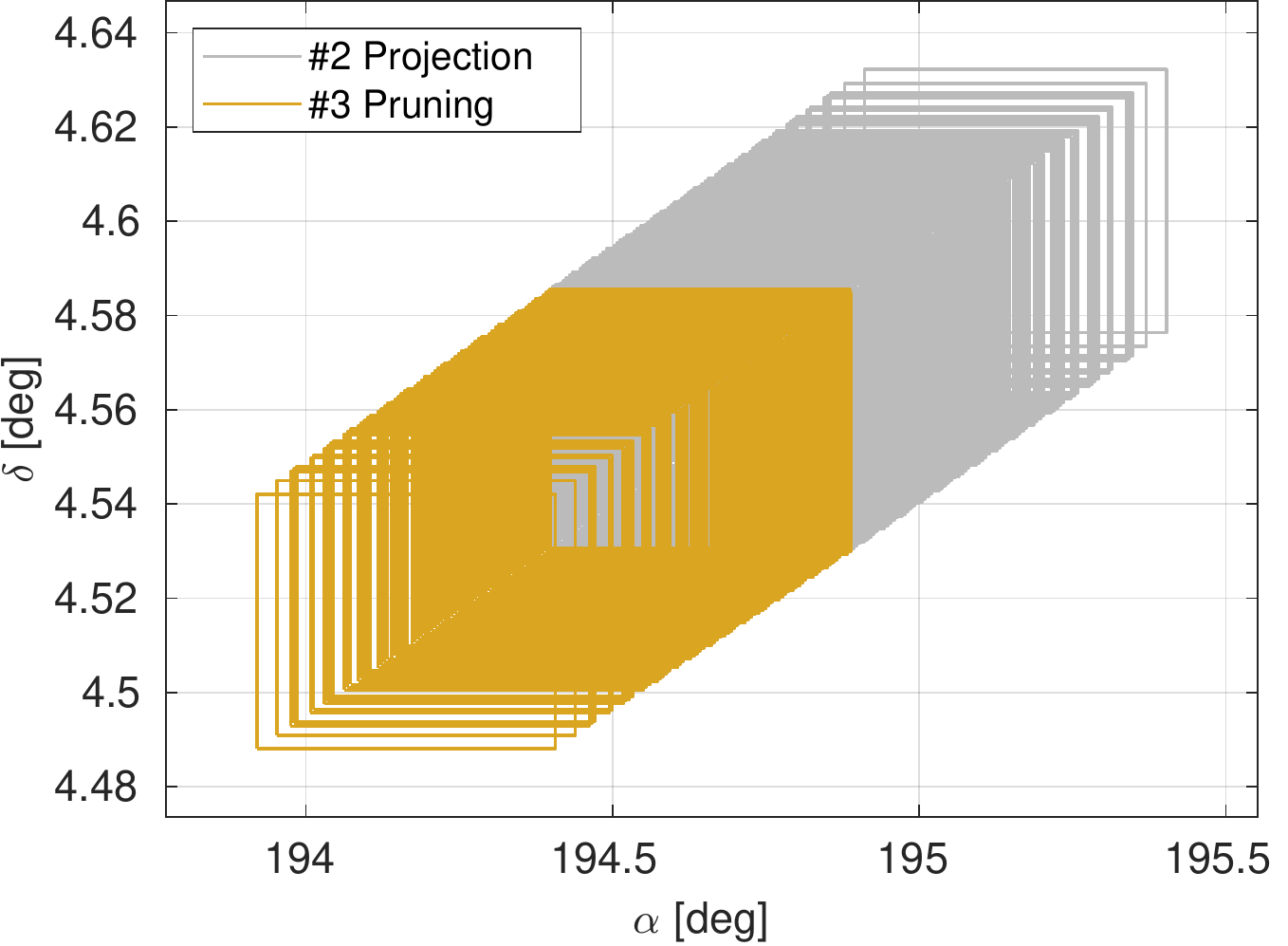}
    \captionsetup{justification=centering}
    \caption{Estimated bounds of projected and retained domains at \SI{100.366}{\hour} after $t_0$ for scenario H}
    \label{fig:bounds_true_out}
\end{figure}

Two examples of pruning are depicted in \cref{fig:bounds_true_no_out,fig:bounds_true_out} for scenarios F and H respectively. These correspond to the epochs with the highest ratio of projected over retained domains, i.e. when the largest number of domains is dropped and the uncertainty is reduced the most. As expected, for scenario H this occurs at the first measurements after the block of outliers which ends a long propagation arc without useful information to reduce the growing uncertainty (since uncorrelated measurements are ignored).

\begin{figure}[!ht]
    \centering
    \includegraphics[width=0.6\textwidth]{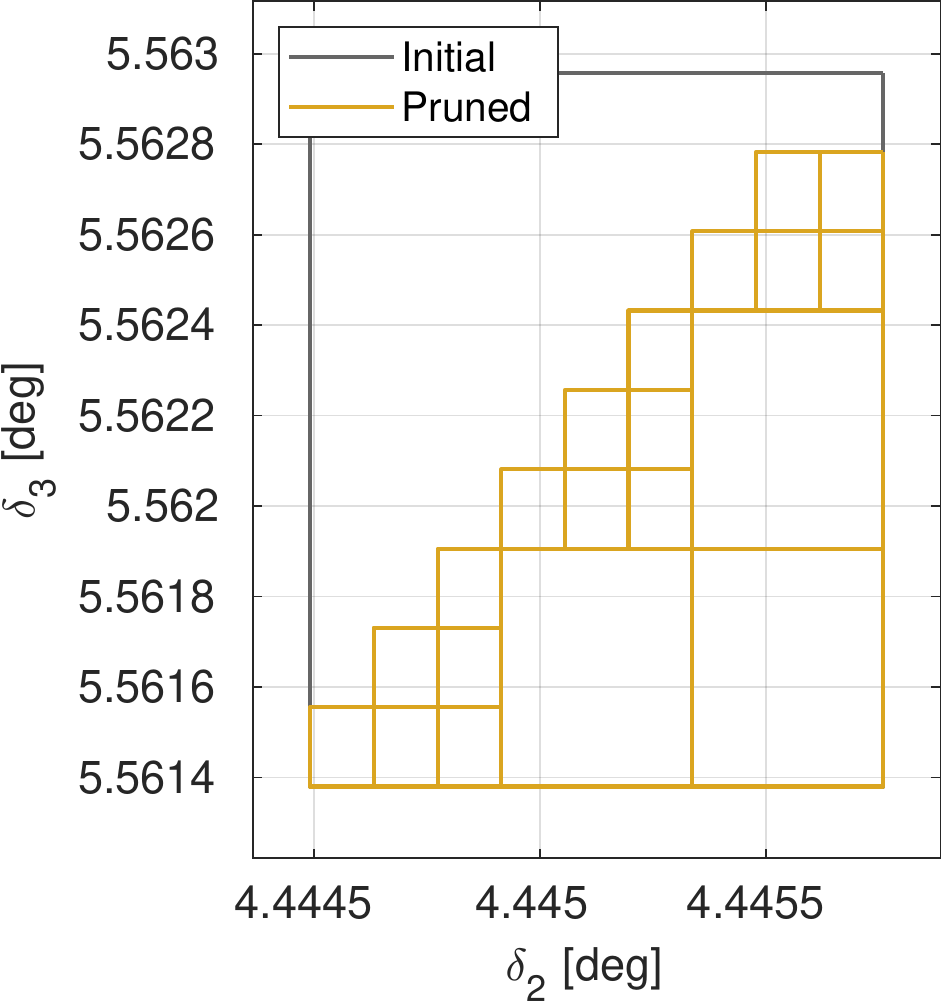}
    \captionsetup{justification=centering}
    \caption{Edges of the base \glsentryshort{iod} solution and retained domains in $(\alpha_2,\delta_2)$ space for scenario F}
    \label{fig:iod_true_no_out}
\end{figure}

\begin{figure}[!ht]
    \centering
    \includegraphics[width=0.6\textwidth]{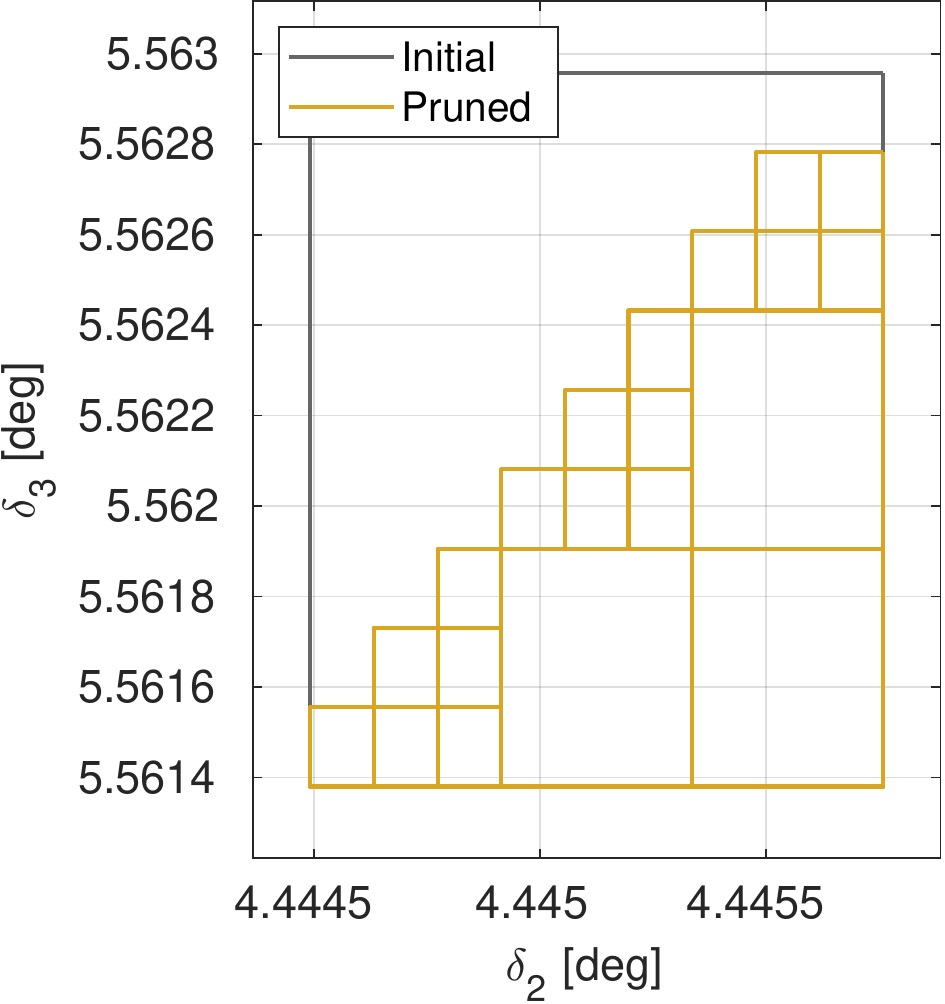}
    \captionsetup{justification=centering}
    \caption{Edges of the base \glsentryshort{iod} solution and retained domains in $(\alpha_2,\delta_2)$ space for scenario H}
    \label{fig:iod_true_out}
\end{figure}

The uncertainty on the \glsentryshort{iod} solution for the two scenarios is depicted in \cref{fig:iod_true_no_out,fig:iod_true_out} in which a projection of the initial (gray box) and retained domains (golden boxes) onto the $(\delta_2,\delta_3)$ space is shown. A comparison of the area covered by the golden boxes demonstrates that for both scenarios the correct subregion of the initial domain that contains the true measurement is identified. The last corresponds to the lower-right portion of the initial rectangular domain, i.e. to above-average values of $\delta_2$ and below-average values of $\delta_3$.

Similarly to \cref{sec:synthetic_measurements}, \glsentryshort{od} is performed in the four scenarios to assess the performance of the \glsentryshort{ls} and \glsentryshort{lsar} estimators while processing both raw and retained observations. The true state is not available in this case, and the best estimate available from the \glsentryshort{cnes}'s catalog is used as ground truth. The last was obtained by processing measurements spanning several years of observations, and it is continuously updated as soon as new data is available using a very \glsentrylong{hf} dynamical model to propagate the best estimate to each measurement epoch. For the shake of limiting the overall runtime, this paper considers only a very limited number of observations, and no prior information is assumed except that all measurements are correlated to the same \glsentryshort{so}. The \glsentrylong{hf} dynamical model was also simplified for the same reason. As a consequence, errors in the order of \SI{10}{\km} and \SI{0.001}{\km\per\s} over a five days time span are considered acceptable in this case. The aforementioned reference solution is thus reported in \cref{tab:od_true_true} together with the initial guesses provided to the batch estimators. The lasts coincide with the raw \glsentryshort{iod} solution for scenarios E and G (since no pruning is performed) and to the domain that is characterized by the smallest residuals among the retained ones for scenarios F and H.

\begin{table}[!ht]
    \centering
    \captionsetup{justification=centering}
    \caption{True state and initial guesses for scenarios E to H}
    \begin{adjustwidth}{-2.5cm}{-2.5cm}
    \centering
    \small
    \sisetup{round-mode=places,round-precision=5}
    \setlength\tabcolsep{2pt}
    \begin{tabular}{c *6{S[scientific-notation=true,table-format=-1.5e-1]}}
    \toprule
    & {$x$, \si{\km}} & {$y$, \si{\km}} & {$z$, \si{\km}} & {$v_x$, \si{\km\per\s}} & {$v_y$, \si{\km\per\s}} & {$v_z$, \si{\km\per\s}}\\ \midrule
    True & -21551.1846646302 & 14404.8664520748 & -1082.46255877053 & -3.58040390149113 & -0.736464589895252 & 0.00194379476500789\\ \midrule
    Guess E & -21547.0726142591 & 14404.6561866209 & -1081.990386528 & -3.58085716435414 & -0.735899837761463 & 0.00200491548637004\\
    Guess F & -21551.8866167486 & 14407.4165534004 & -1081.65956618829 & -3.58012190214936 & -0.736320077278301 & 0.00206423491233456\\
    Guess G & -21547.0726142591 & 14404.6561866209 & -1081.990386528 & -3.58085716435414 & -0.735899837761463 & 0.00200491548637004\\
    Guess H & -21551.8866167486 & 14407.4165534004 & -1081.65956618829 & -3.58012190214936 & -0.736320077278301 & 0.00206423491233456\\ \bottomrule
    \end{tabular}
    \end{adjustwidth}
    \label{tab:od_true_true}
\end{table}

\Cref{tab:od_sol_true_no_out} reports the \glsentryshort{ls} and \glsentryshort{lsar} estimation errors and $\pm 3\sigma$ bounds for scenario E. These are qualitatively identical to those of scenario F which are omitted for conciseness. The \glsentryshort{ls} solution is again characterized by tighter bounds on the states, even though the errors with respect to the assumed reference are larger than for the \glsentryshort{lsar} estimate. In this case, three iterations are needed for both estimators to converge.

\begin{table}[!ht]
    \centering
    \captionsetup{justification=centering}
    \caption{\glsentryshort{ls} and \glsentryshort{lsar} estimation errors and $\pm 3\sigma$ bounds for scenario E}
    \begin{adjustwidth}{-2.5cm}{-2.5cm}
    \centering
    \small
    \sisetup{round-mode=places,round-precision=3}
    \setlength\tabcolsep{2pt}
    \begin{tabular}{r *6{S[scientific-notation=true,table-format=1.3e-1]}}
    \toprule
    & {$x$, \si{\km}} & {$y$, \si{\km}} & {$z$, \si{\km}} & {$v_x$, \si{\km\per\s}} & {$v_y$, \si{\km\per\s}} & {$v_z$, \si{\km\per\s}}\\ \midrule
    \glsentryshort{ls} Error & 9.2394503500052 & 3.30574014850232 & 0.116199529730007 & 0.00175312153853274 & 0.000306294804796094 & 0.000136508035645701\\
    $\pm 3\sigma$ bounds & 0.780987256152968 & 0.41443751066775 & 0.138698947819167 & 0.000150296943704622 & 3.3007112216976e-05 & 1.08482713045251e-05\\ \midrule
    \glsentryshort{lsar} Error & 0.602359107738526 & 2.32368489829035 & 0.774951866708533 & 0.000183186873129161 & 7.60796399684518e-05 & 0.000171030451399098\\
    $\pm 3\sigma$ bounds & 0.966984353070529 & 0.513004613089203 & 0.172032715398441 & 0.000186008968325118 & 4.08538660639667e-05 & 1.34418171818632e-05\\ \bottomrule
    \end{tabular}
    \end{adjustwidth}
    \label{tab:od_sol_true_no_out}
\end{table}

The effects of measurements outliers among the processed set is then shown in \cref{tab:od_sol_true_out_no_prune} which summarizes the results obtained for scenario G. In this case, five and three iterations are needed by the \glsentryshort{ls} and \glsentryshort{lsar} algorithms to converge. As for scenario C, the \glsentryshort{ls} estimate exhibits large errors in both positions ($\approx\SI{550}{\km}$) and velocity ($\approx\SI{0.1}{\km\per\s}$) thus supporting the needs for a preprocessing algorithm capable of detecting and isolating uncorrelated observations. On the contrary, the \glsentryshort{lsar} scheme shows an excellent performance that matches the one of scenario H reported in \cref{tab:od_sol_true_out_prune}. At the same time, the presence of uncorrelated measurements among the processed set does not influence the estimated variance as can be verified by comparing the second and fourth rows in \cref{tab:od_sol_true_no_out,tab:od_sol_true_out_no_prune} respectively.

\begin{table}[!ht]
    \centering
    \captionsetup{justification=centering}
    \caption{\glsentryshort{ls} and \glsentryshort{lsar} estimation errors and $\pm 3\sigma$ bounds for scenario G}
    \begin{adjustwidth}{-2.5cm}{-2.5cm}
    \centering
    \small
    \sisetup{round-mode=places,round-precision=3}
    \setlength\tabcolsep{2pt}
    \begin{tabular}{r *6{S[scientific-notation=true,table-format=1.3e-1]}}
    \toprule
    & {$x$, \si{\km}} & {$y$, \si{\km}} & {$z$, \si{\km}} & {$v_x$, \si{\km\per\s}} & {$v_y$, \si{\km\per\s}} & {$v_z$, \si{\km\per\s}}\\ \midrule
    \glsentryshort{ls} Error & 464.561341965927 & 289.620372040554 & 31.6615004965373 & 0.0937668132322371 & 0.0205339834402047 & 0.00156537307514134\\
    $\pm 3\sigma$ bounds & 0.767690130097545 & 0.409522817594827 & 0.143979355581707 & 0.0001455951781223 & 3.1723125929194e-05 & 1.06989268261277e-05\\ \midrule
    \glsentryshort{lsar} Error & 1.25560162653959 & 2.69372598359898 & 0.819188049874167 & 0.000288182607168274 & 2.51686165960147e-05 & 0.000180207964448427\\
    $\pm 3\sigma$ bounds & 1.00660847682005 & 0.540640712258467 & 0.176036854759609 & 0.000194414857751158 & 4.26913791695419e-05 & 1.37102004366034e-05\\ \bottomrule
    \end{tabular}
    \end{adjustwidth}
    \label{tab:od_sol_true_out_no_prune}
\end{table}

\begin{table}[!ht]
    \centering
    \captionsetup{justification=centering}
    \caption{\glsentryshort{ls} and \glsentryshort{lsar} estimation errors and $\pm 3\sigma$ bounds for scenario H}
    \begin{adjustwidth}{-2.5cm}{-2.5cm}
    \centering
    \small
    \sisetup{round-mode=places,round-precision=3}
    \setlength\tabcolsep{2pt}
    \begin{tabular}{r *6{S[scientific-notation=true,table-format=1.3e-1]}}
    \toprule
    & {$x$, \si{\km}} & {$y$, \si{\km}} & {$z$, \si{\km}} & {$v_x$, \si{\km\per\s}} & {$v_y$, \si{\km\per\s}} & {$v_z$, \si{\km\per\s}}\\ \midrule
    \glsentryshort{ls} Error & 3.29604513153692 & 0.32672932171563 & 0.505709718794199 & 0.000549673851591107 & 4.28823231984432e-05 & 0.000156495686685633\\
    $\pm 3\sigma$ bounds & 0.923162975723824 & 0.512041529535699 & 0.14964674270472 & 0.000180287872758937 & 3.95567781780505e-05 & 1.153318460807e-05\\ \midrule
    \glsentryshort{lsar} Error & 1.1518140723987 & 2.63520422168635 & 0.815017315887631 & 0.000268448717855625 & 2.93810689749442e-05 & 0.000179925964703669\\
    $\pm 3\sigma$ bounds & 1.14455295583921 & 0.634817164172961 & 0.18559548708244 & 0.000223489960097547 & 4.90322609949816e-05 & 1.42958552451257e-05\\ \bottomrule
    \end{tabular}
    \end{adjustwidth}
    \label{tab:od_sol_true_out_prune}
\end{table}

\begin{table}[!ht]
    \centering
    \captionsetup{justification=centering}
    \caption{Runtime and number of iterations with real measurements}
    \begin{adjustwidth}{-2.5cm}{-2.5cm}
    \centering
    \small
    \sisetup{round-mode=places,round-precision=3}
    \setlength\tabcolsep{2pt}
    \begin{tabular}{c *3{S[scientific-notation=false,table-format=3.3]} cc}
    \toprule
    & \multicolumn{3}{c}{$t$, \si{\s}} & \multicolumn{2}{c}{\# iterations}\\ \midrule
    Case \# & Pruning & LS & LSAR & LS & LSAR\\ \midrule
    E & {-} & 510.904 & 200.603 & 3 & 3\\
    F & 49.653 & 453.881 & 198.409 & 3 & 3\\
    G & {-} & 537.406 & 192.47 & 5 & 3\\
    H & 14.897 & 464.756 & 140.941 & 5 & 3\\ \bottomrule
    \end{tabular}
    \end{adjustwidth}
    \label{tab:od_runtime_tarot}
\end{table}

Looking at \cref{tab:od_runtime_tarot}, no reduction in the number of iterations is observed when processing correlated measurements only instead of the entire available set. The inclusion of the pruning step is however paramount to restore the accuracy of the \glsentryshort{ls} estimate as demonstrated in \cref{tab:od_sol_true_out_prune}. In this case a consistent estimate is in fact obtained and the \glsentryshort{ls} estimator converges to a solution closer to the true one than that of \cref{tab:od_sol_true_no_out}.

\section{Conclusions}\label{sec:conclusions}

A pruning algorithm for the refinement of an \glsentryshort{iod} solution and the detection of measurement outliers was presented. This scheme is broken down into four main steps, carried out sequentially for each available observation: propagation of the state estimate, projection onto the observables space, domains pruning based on the matching between predicted and actual measurements, and domains merging to minimize the number of trajectories to be further propagated. The uncertainty on the \glsentryshort{iod} solution is represented as a Taylor polynomial and \glsentryshort{da} techniques are leveraged to efficiently propagate the uncertainty over time. A novel splitting algorithm is also employed for the automatic control of the error committed by the truncated polynomials. To further reduce its computational effort, the prediction step is based on a bifidelity dynamical model that couples a \glsentryshort{da}-aware \glsentryshort{sgp} propagator with an \glsentrylong{hf} numerical integrator. If the first is used to map the initial uncertainty on the state to the next observation epoch, the last is employed to compute an accurate reference trajectory on which said uncertainty is centered and to include the effects of stochastic accelerations in the predicted estimate. After projecting the prior distribution onto the observables space, polynomial bounding techniques are used to estimate the range of variation of each domain and prune the ones that do not overlap with observed data at the same epoch, thus reducing the volume of the propagated uncertainty and tightening the bounds on the initial \glsentryshort{iod} solution. Measurements outliers are also detected in this step if no intersection between predicted and actual measurements is found. This algorithm is then employed as a preprocessing scheme to filter the information used by subsequent batch \glsentryshort{od} tools. Its performance is demonstrated for a \glsentryshort{gto} object using synthetic as well as real observation data from the \glsentryshort{tarot} network. It is verified to greatly improve the accuracy of \glsentryshort{ls} estimates, notoriously affected by the presence of outliers among the processed observations, and to reduce the runtime of \glsentryshort{lsar} estimators, which converge in less iterations if processing correlated measurements only instead of the entire available set. The successful application of this method to real data in a challenging scenario, namely a \glsentryshort{heo} object characterized by poor observability and subject to a strongly nonlinear dynamics, suggests that it is applicable to the full spectrum of orbit regimes.

\section*{Acknowledgments}

This work is co-funded by the \glsentryshort{cnes} through A. Foss\`a PhD program, and made use of the \glsentryshort{cnes} orbital propagation tools, including the \glsentryshort{pace} library. The authors gratefully acknowledge Dr. Laura Pirovano for her support during the implementation and validation of the \glsentryshort{iod} algorithm.

\printbibliography

\end{document}